\documentstyle[12pt,twoside]{article}

\font\teneufm=eufm10 scaled \magstep1
\font\seveneufm=eufm7 scaled \magstep1
\font\fiveeufm=eufm5  scaled \magstep1
\newfam\eufmfam
\textfont\eufmfam=\teneufm
\scriptfont\eufmfam=\seveneufm
\scriptscriptfont\eufmfam=\fiveeufm
\def\frak#1{{\fam\eufmfam\relax#1}}

\newfam\msbfam
\font\tenmsb=msbm10 scaled \magstep1  \textfont\msbfam=\tenmsb
\font\sevenmsb=msbm7 scaled \magstep1 \scriptfont\msbfam=\sevenmsb
\font\fivemsb=msbm5 scaled \magstep1  \scriptscriptfont\msbfam=\fivemsb
\def\Bbb{\fam\msbfam \tenmsb}

\def\RR{{\Bbb R}}
\def\CC{{\Bbb C}}
\def\QQ{{\Bbb Q}}
\def\NN{{\Bbb N}}
\def\ZZ{{\Bbb Z}}

\def\PP{{\Bbb P}}

\def\ra{\rightarrow}

 \def\HollowBoxx #1#2#3{{\dimen0=#1 \advance\dimen0 by -#2
       \dimen1=#1 \advance\dimen1 by #3
        \vrule height 0pt depth #3 width #2
       \hskip -#3
       \vrule height #1 depth #3 width #3}}
 \def\LeftContraction{\mathord{\kern1.45pt \HollowBoxx{6pt}{3.5pt}{.4pt}}\,}

 \def\HollowBox #1#2#3{{\dimen0=#1 \advance\dimen0 by -#3
       \dimen1=#1 \advance\dimen1 by #3
        \vrule height #1 depth #3 width #3
        \vrule height 0pt depth #3 width #2
        \hskip -#3}}
 \def\RightContraction{\mathord{\, \HollowBox{6pt}{3.1pt}{.4pt}} \kern1.6pt}

\def\qed{{\hfill $\Box$}}
\newtheorem{theorem}{THEOREM}[section]

\newtheorem{lemma}[theorem]{Lemma}
\newtheorem{example}[theorem]{Example}
\newtheorem{remark}[theorem]{Remark}
\newtheorem{proposition}[theorem]{Proposition}

\newtheorem{definition}[theorem]{Definition}

\begin{document}

\begin{center}
{\Large \bf Effective Actions of the Unitary Group
\medskip\\
 on Complex Manifolds}
\medskip \\
\normalsize A. V. Isaev and N. G. Kruzhilin
\end{center}

\begin{quotation} \small \sl We classify all connected $n$-dimensional
  complex manifolds admitting an effective action of the unitary
  group $U_n$ by biholomorphic transformations. One
consequence of this
classification is a
  characterization of  $\CC^n$ by its
  automorphism group.
\end{quotation}

\pagestyle{myheadings}
\markboth{A. V. Isaev and N. G. Kruzhilin}{Unitary Group Actions}

\setcounter{section}{-1}

\section{Introduction}
\setcounter{equation}{0}

We are interested in classifying
 all connected complex manifolds $M$ of
dimension $n\ge 2$ admitting
effective actions of the unitary group $U_n$ by biholomorphic
transformations.

One motivation for our study was the
following question that we learned from S. Krantz: assume that  the group
$\hbox{Aut}(M)$ of all biholomorphic automorphisms of $M$ and the group
$\hbox{Aut}(\CC^n)$ of all biholomorphic automorphisms of $\CC^n$ are isomorphic
as topological groups equipped with the compact-open
 topology; does it
imply that $M$ is biholomorphically equivalent to $\CC^n$? The group
$\hbox{Aut}(\CC^n)$ is very large (see, e.g., \cite{AL}), and it is not
that clear from 
the start what automorphisms of $\CC^n$ one can  use to
approach the problem. The isomorphism
between $\hbox{Aut}(M)$ and $\hbox{Aut}(\CC^n)$ induces a continuous
effective action on $M$ of any subgroup $G\subset\hbox{Aut}(\CC^n)$. If $G$ is
a Lie group, then this action is in fact real-analytic. We consider $G=U_n$
which, as it turns out, results in a very short list of manifolds that can
occur.

In Section 1 we find all possible dimensions of orbits  of a
 $U_n$-action on $M$.  It turns out (see Proposition \ref{dim}) that an
orbit is  either  a point (hence $U_n$ has a fixed point in
$M$), or a real hypersurface in $M$, or a complex hypersurface in $M$,
or the whole of $M$ (in which case $M$ is homogeneous).

Manifolds admitting an action with  fixed point were found in
\cite{K} (see Remark \ref{kaup}).

In Section 2 we classify
manifolds with  a $U_n$-action
such that all orbits are real hypersurfaces. We show that
such a  manifold  is  either a spherical layer in $\CC^n$, or a Hopf
manifold, or the quotient  of one of these manifolds
 by the action of a discrete subgroup of the center of
$U_n$ (Theorem \ref{finalstep}).

In Section 3 we consider the situation when every orbit is  a
real or a complex hypersurface in $M$ and show that there can exist  at most
 two orbits that are complex hypersurfaces.
Moreover,  such  orbits turn out to be
 biholomorphically equivalent to $\CC\PP^{n-1}$ and can only
arise either as a result of blowing up
 $\CC^n$ or a ball in $\CC^n$
 at the origin, or adding the hyperplane
$\infty\in\CC\PP^n$ to the exterior of a ball in $\CC^n$,
or blowing up $\CC\PP^n$ at
one point,  or taking the quotient of one
of these examples by the action  of a discrete subgroup of the center
of $U_n$ (Theorem \ref{complh1}).

Finally, in Section 4 we consider the homogeneous case.
In this case the manifold in question must   be
equivalent to the quotient of a Hopf manifold
 by the action of a discrete central subgroup
  (Theorem
\ref{hopfclass}).

Thus, Remark \ref{kaup}, Theorem \ref{finalstep}, Theorem
\ref{complh1} and Theorem \ref{hopfclass} provide a complete list of
connected manifolds of dimension $n\ge 2$ admitting effective actions
of $U_n$ by biholomorphic transformations. An easy consequence of this
classification is the following characterization of $\CC^n$ by its
automorphism group that we obtain in Section 5:
\smallskip\\

\noindent {\bf THEOREM \ref{char}} {\sl Let $M$ be a connected complex manifold of dimension
  $n$. Assume that $\hbox{Aut}(M)$ and
  $\hbox{Aut}(\CC^n)$ are isomorphic as topological groups. Then $M$
  is biholomorphically equivalent to $\CC^n$.}
\smallskip\\

We acknowledge that this work started while the second author
was  visiting Centre for Mathematics and its Applications,
Australian National University.

\section{Dimensions of Orbits}
\setcounter{equation}{0}

In this section we obtain the following result, which is similar to
Satz 1.2 in \cite{K}.

\begin{proposition}\label{dim} \sl Let $M$ be a connected complex manifold of
  dimension $n\ge 2$ endowed with an effective action of   $U_n$
by biholomorphic transformations. Let $p\in M$ and
let $O(p)$ be  the $U_n$-orbit of $p$. Then $O(p)$
  is either

\noindent (i) the whole of $M$ (hence $M$ is compact), or

\noindent (ii) a single point, or

\noindent (iii) a complex compact hypersurface in $M$, or

\noindent (iv) a real compact hypersurface in $M$.

\end{proposition}

\noindent {\bf Proof:} For $p\in M$
  let  $I_p$ be  the isotropy subgroup of $U_n$ at $p$, i.e.,
  $I_p:=\{g\in U_n: gp=p\}$. We denote by $\Psi$
the continuous homomorphism of
$U_n$ into $\hbox{Aut}(M)$ (the group of biholomorphic automorphisms of
$M$) induced by the action of $U_n$ on $M$.
Let $L_p:=\{d_p(\Psi(g)): g\in I_p\}$ be the
  linear isotropy subgroup, where $d_pf$ is the differential of a
  map $f$ at $p$.
 Clearly,
  $L_p$ is a compact subgroup of $GL(T_p(M),\CC)$. Since the action of
  $U_n$ is effective, $L_p$ is isomorphic to
  $I_p$. Let $V\subset T_p(M)$ be the tangent space to $O(p)$ at $p$. Clearly, $V$ is
   $L_p$-invariant. We assume now that $O(p)\ne M$ (and therefore $V\ne
  T_p(M)$)
and consider the following three cases.
\smallskip\\

{\bf Case 1.}  $d:=\hbox{dim}_{\CC}(V+iV)<n$.
\smallskip\\

Since $L_p$ is compact, one can consider coordinates on $T_p(M)$  such
 that $L_p\subset U_n$. Further, the action of $L_p$ on $T_p(M)$
is completely reducible and the subspace $V+iV$ is invariant  under this
action. Hence  $L_p$ can in fact be embedded in $U_d\times
U_{n-d}$. Since $\hbox{dim}\,O(p)\le 2d$, it follows that
$$
n^2\le d^2+(n-d)^2+2d,
$$
and therefore  either $d=0$ or $d=n-1$. If $d=0$, then we obtain
  (ii). If $d=n-1$, then  the above relation is in fact the equality
$\hbox{dim}\,O(p)=2d=2n-2$, and therefore $iV=V$,
which yields (iii).
\smallskip\\

{\bf Case 2.}  $T_p(M)=V+iV$ and $r:=\hbox{dim}_{\CC}(V\cap iV)>0$.
\smallskip\\

As above, $L_p$ can be embedded in $U_r\times U_{n-r}$ (clearly, we have
$r<n$).  Moreover,
 $V\cap iV\ne V$ and since $L_p$ preserves $V$, it follows that
$\hbox{dim}\,L_p<r^2+(n-r)^2$. We have $\hbox{dim}\,O(p)\le 2n-1$, and
therefore
$$
n^2<r^2+(n-r)^2+2n-1,
$$
which shows  that $\hbox{dim}\,O(p)=2n-1$. This yields (iv).
\smallskip\\

{\bf Case 3.}  $T_p(M)=V\oplus iV$.
\smallskip\\

In this case  $\hbox{dim}\, V=n$ and
$L_p$ can be embedded in the real orthogonal group $O_n(\RR)$,
and therefore
$$
\hbox{dim}\,L_p+\hbox{dim}\, O(p)\le \frac{n(n-1)}{2}+n<n^2,
$$
which is a contradiction.

The proof of the proposition  is complete.\qed
\smallskip\\

\begin{remark}\label{kaup}\rm  It is shown in \cite{K} (see Folgerung
  1.10 there)
  that if $U_n$ has a fixed point in $M$, then $M$ is
  biholomorphically equivalent to
  either

\noindent (i) the unit ball $B^n\subset\CC^n$, or

\noindent (ii) $\CC^n$, or

\noindent (iii) $\CC\PP^n$.

\noindent The biholomorphic equivalence $f$ can be
  chosen to be an isomorphism of $U_n$-spaces, more precisely,
$$
f(gq)=\gamma(g)f(q),
$$
where either $\gamma(g)=g$ or $\gamma(g)=\overline{g}$ for all $g\in U_n$ and $q\in M$
(here $B^n$, $\CC^n$ and $\CC\PP^n$ are
considered with the standard actions of $U_n$).
\end{remark}

\section{The Case of Real Hypersurface Orbits}
\setcounter{equation}{0}

We shall now consider  orbits in
$M$ that are  real hypersurfaces. We require the following algebraic result.

\begin{lemma}\label{un} \sl Let $G$ be a connected closed subgroup of
  $U_n$ of dimension $(n-1)^2$, $n\ge 2$. Then either $G$ contains
  the center of $U_n$, or $G$ is conjugate in $U_n$ to
   the subgroup of all matrices
\begin{equation}
\left(\begin{array}{cc}
\alpha & 0\\
0& \beta
\end{array}\right),\label{spec}
\end{equation}
where $\alpha \in U_1$ and $\beta\in SU_{n-1}$, or
for some  $k_1,k_2\in\ZZ$, $(k_1,k_2)=1$, $k_2\ne 0$, it is conjugate
 to the subgroup
$H_{k_1,k_2}$ of all matrices
\begin{equation}
\left(\begin{array}{cc}
a & 0\\
0 & B
\end{array}\right),\label{mat}
\end{equation}
where $B\in U_{n-1}$ and $a\in
(\det B)^{\frac{k_1}{k_2}}:=\exp(k_1/k_2\, \hbox{Ln}\,(\det B))$.
\end{lemma}

\noindent {\bf Proof:} Since $G$ is compact, it is completely
reducible, i.e., $\CC^n$ splits into a sum of $G$-invariant
pairwise orthogonal complex subspaces, $\CC^n=V_1\oplus\dots\oplus V_m$,
such that the restriction $G_j$ of $G$ to each  $V_j$ is irreducible. Let
$n_j:=\hbox{dim}_{\CC}V_j$ (hence $n_1+\dots+n_m=n$) and let
$U_{n_j}$ be the group of unitary
transformations of $V_j$. Clearly, $G_j\subset U_{n_j}$, and therefore
$\hbox{dim}\,G\le n_1^2+\dots+n_m^2$.
On the other hand  $\hbox{dim}\,G=(n-1)^2$, which shows that
$m\le 2$.

Let $m=2$. Then there exists a unitary change of coordinates
$\CC^n$ such that in the new variables   elements of $G$ are of  the
 form
\begin{equation}
\left(\begin{array}{cc}
a & 0\\
0& B
\end{array}\right),\label{mat1}
\end{equation}
where $a \in U_1$ and $B \in U_{n-1}$. We note that the scalars  $a$
and the matrices $B$  in (\ref{mat1}) corresponding to the elements of $G$
form
compact connected subgroups of $U_1$ and $U_{n-1}$, respectively; we
shall denote  them by
$G_1$ and $G_2$ as above.

If $\hbox{dim}\,G_1=0$, then $G_1=\{1\}$, and
therefore $G_2=U_{n-1}$. Thus we get
the form (\ref{mat}) with $k_1=0$.

Assume  that $\hbox{dim}\,G_1=1$, i.e., $G_1=U_1$. Then
$(n-1)^2-1\le\hbox{dim}\,G_2\le (n-1)^2$.  Let
$\hbox{dim}\,G_2=(n-1)^2-1$ first. The only connected subgroup of $U_{n-1}$
of dimension $(n-1)^2-1$ is $SU_{n-1}$. Hence $G$ is conjugate to
the subgroup
of matrices of the form (\ref{spec}).
 Now let  $\hbox{dim}\,G_2=(n-1)^2$, i.e., $G_2=U_{n-1}$.
Consider the Lie algebra ${\frak g}$ of
  $G$. It consists of matrices of the following form:
\begin{equation}
\left(\begin{array}{cc}
l(b)& 0\\
0& b
\end{array}\right),\label{mat2}
\end{equation}
where $b$ is an arbitrary matrix in
 ${\frak u}_{n-1}$ and $l(b)\not\equiv 0$
is a linear
function  of the matrix elements of $b$ ranging  in $i\RR$. Clearly,
$l(b)$ must vanish on the commutant of ${\frak u}_{n-1}$,
 which is
${\frak {su}}_{n-1}$. Hence matrices (\ref{mat2}) form  a
Lie algebra if and only if $l(b)=c\cdot\hbox{trace}\,b$, where
$c\in\RR\setminus\{0\}$. Such an algebra can be the Lie algebra of
a subgroup of $U_1\times U_{n-1}$  only if
$c\in\QQ\setminus\{0\}$.
Hence $G$ is conjugate to the group of matrices
(\ref{mat}) with some $k_1,k_2\in\ZZ$, $k_2\ne 0$, and one can always
assume that $(k_1,k_2)=1$.

Now let $m=1$. We  shall proceed as in the proof of Lemma 2.1 in \cite{IK}.
Let ${\frak g}\subset{\frak {u}}_n\subset{\frak {gl}}_n$ be the Lie algebra of $G$ and
${\frak g}^{\CC}:={\frak g}+i{\frak g}\subset{\frak {gl}}_n$ its complexification.
Then
 ${\frak g}^{\CC}$ acts irreducibly on $\CC^n$ and
by a theorem of \'E. Cartan (see, e.g., \cite{GG}), ${\frak g}^{\CC}$ is
either semisimple or  the direct sum of a semisimple ideal ${\frak
h}$ and the center of ${\frak {gl}}_n$ (which is isomorphic to $\CC$).
Clearly, the action of the ideal ${\frak h}$ on $\CC^n$ must be  irreducible.

Assume first that ${\frak g}^{\CC}$ is semisimple, and let ${\frak
g}^{\CC}={\frak g}_1\oplus\dots\oplus{\frak g}_k$ be its decomposition
into the direct sum of simple ideals. Then  (see,
e.g., \cite{GG}) the irreducible $n$-dimensional representation of
${\frak g}^{\CC}$ given by the embedding of ${\frak g}^{\CC}$ in ${\frak {gl}}_n$ is the
tensor product of some irreducible faithful representations of the ${\frak
g}_j$. Let $n_j$ be the dimension of the corresponding representation  of ${\frak
g}_j$, $j=1,\dots,k$.
 Then $n_j\ge 2$, $\hbox{dim}_{\CC}\,{\frak g}_j\le n_j^2-1$,
and $n=n_1\cdot\dots\cdot n_k$.
The following observation is simple.

\begin{quote}
{\bf Claim:} {\sl If $n=n_1\cdot\dots\cdot n_k$, $k\ge 2$, $n_j\ge 2$
for $j=1,\dots,k$, then $\sum_{j=1}^k n_j^2\le n^2-2n$.}
\end{quote}

Since $\hbox{dim}_{\CC}\,{\frak g}^{\CC}=(n-1)^2$, it follows from the
above claim that $k=1$, i.e., ${\frak g}^{\CC}$ is
simple. The minimal dimensions of irreducible faithful
representations of complex simple Lie algebras are well-known (see,
e.g., \cite{VO}). In the table below $V$ denotes
 representations of minimal dimension.

\begin{center}
\begin{tabular}{|l|c|c|}
\hline
\multicolumn{1}{|c|}{${\frak g}$}&
\multicolumn{1}{c|}{ $\hbox{dim}\,V$}&
\multicolumn{1}{c|}{$\hbox{dim}\,{\frak g}$}
\\ \hline
${\frak {sl}}_k$\,\,$k\ge 2$ & $k$ & $k^2-1$
\\ \hline
${\frak o}_k$\,\, $k\ge 7$&  $k$ & $\frac{k(k-1)}{2}$
\\ \hline
${\frak {sp}}_{2k}$\,\,$k\ge 2$ & $2k$ & $2k^2+k$
\\ \hline
${\frak e}_6$ & 27 & 78
\\ \hline
${\frak e}_7$ & 56 & 133
\\ \hline
${\frak e}_8$ & 248 & 248
\\ \hline
${\frak f}_4$ & 26 & 52
\\ \hline
${\frak g}_2$ & 7 & 14
\\ \hline
\end{tabular}
\end{center}

Since $\hbox{dim}_{\CC}\,{\frak g}^{\CC}=(n-1)^2$, it follows that
none  of the above possibilities realize. Hence  ${\frak g}^{\CC}$ contains the center of ${\frak {gl}}_n$,
and therefore ${\frak g}$ contains the center of ${\frak u}_n$. Thus
 $G$ contains the center of $U_n$.

The proof of the lemma is complete.\qed
\smallskip\\

We can now prove the following proposition.

\begin{proposition}\label{realh} \sl Let $M$ be a
  complex manifold of dimension $n\ge 2$ endowed with an effective
action of   $U_n$ by
   biholomorphic transformations. Let
  $p\in M$ and let the orbit $O(p)$ be a real hypersurface in $M$. Then
the isotropy subgroup
  $I_p$ is isomorphic to $U_{n-1}$.
\end{proposition}

\noindent {\bf Proof:} Since $O(p)$ is a real hypersurface in $M$, it
  arises in Case 2 in the proof of Proposition \ref{dim}. We shall
use the notation from that proof. Let $W$ be
  the orthogonal complement to $V\cap iV$ in $T_p(M)$. Clearly,
  $\hbox{dim}_{\CC}\,V\cap iV=n-1$ and
  $\hbox{dim}_{\CC}\,W=1$. The group $L_p$ is a subgroup of
  $U_n$ and preserves $V$, $V\cap iV$,  and
  $W$;  hence it preserves  the line $W\cap
  V$. Therefore, it can act only  as
  $\pm\hbox{id}$ on $W$.
Since $\hbox{dim}\,L_p=(n-1)^2$, the identity component $L_p^c$ of $L_p$
must in fact                     be
  the group of all unitary transformations  preserving $V\cap iV$
  and acting trivially on $W$. Thus,
 $L_p^c$ is isomorphic to $U_{n-1}$ and acts
  transitively on directions in $V\cap iV$. Hence $O(p)$
is either Levi-flat or
  strongly pseudoconvex.

We claim that $O(p)$ cannot be Levi-flat. For assume  that $O(p)$ is Levi-flat.
Then it is foliated by complex
hypersurfaces in $M$. Let ${\frak m}$
be the Lie algebra of all holomorphic vector
fields on $O(p)$ corresponding to the automorphisms of $O(p)$
generated by the action of $U_n$. Clearly, ${\frak m}$ is isomorphic
to ${\frak u}_n$. For $q\in O(p)$ we
denote by $M_q$ the leaf of the foliation passing
through $q$ and consider the subspace ${\frak l}_q\subset{\frak m}$ of
all vector fields tangent to $M_q$ at $q$. Since vector fields in
${\frak l}_q$ remain tangent to $M_q$ at each point in $M_q$, ${\frak
  l}_q$ is in fact a Lie subalgebra of ${\frak m}$.
Clearly, $\hbox{dim}\,{\frak
  l}_q=n^2-1$, and therefore ${\frak l}_q$ is isomorphic to ${\frak
  {su}}_n$. Since there exists only one way to embed ${\frak {su}}_n$
in ${\frak u}_n$, we obtain that the action of $SU_n\subset U_n$
preserves each leaf $M_q$ for $q\in O(p)$. Hence each leaf $M_q$ is a
union of $SU_n$-orbits.
But such an orbit must be open in $M_q$, and therefore
the action of  $SU_n$ is  transitive   on each  $M_q$.

Let $\tilde I_q$ be  the isotropy subgroup of $q$ in
  $SU_n$. Clearly, $\hbox{dim}\,\tilde I_q=(n-1)^2$.
It now follows from Lemma \ref{un} that $\tilde I_q^c$, the connected
identity component of
   $\tilde I_q$, is conjugate in $U_n$ to
the  subgroup $H_{k_1,k_2}$ (see (\ref{mat})) with  $k_1=-k_2=1$.
Hence $\tilde I_q$ contains the
center of $SU_n$.
The elements of the center act trivially on
$SU/\tilde I_q$ (which is equivariantly diffeomorphic to $M_q$). Thus,  the
central elements of $SU_n$ act trivially on each $M_q$, and therefore
on $O(p)$. Consequently, the action of $U_n$ on the real hypersurface $O(p)$,
and therefore on $M$, is not effective,
which is a contradiction showing that  $M$ is strongly pseudoconvex.

Hence $L_p$ can only act identically on $W$. Thus, $L_p$ is
isomorphic to $U_{n-1}$ and so is $I_p$.

The proof is complete.\qed
\smallskip\\

We now classify real hypersurface orbits up to equivariant diffeomorphisms.

\begin{proposition}\label{realh2} \sl Let $M$ be a
  complex manifold of dimension $n\ge 2$ endowed with  an effective
action of  $U_n$
   by biholomorphic transformations. Let
  $p\in M$ and assume  that the orbit
 $O(p)$ is a real hypersurface in $M$. Then
  $O(p)$ is
  isomorphic as a homogeneous space to a lense manifold
${\cal L}^{2n-1}_m:=S^{2n-1}/\ZZ_m$ obtained by identifying each point $x\in S^{2n-1}$ with
  $e^{\frac{2\pi i}{m}}x$, where $m=|nk+1|$, $k\in\ZZ$ 
(here ${\cal L}^{2n-1}_m$ is considered with the standard action of
$U_n/\ZZ_m$).
\end{proposition}

\noindent {\bf Proof:} By Proposition
  \ref{realh}, $I_p$ is isomorphic to $U_{n-1}$.
Hence it follows from Lemma \ref{un} that $I_p$ either contains the
center of $U_n$ or is conjugate to some  group
$H_{k_1,k_2}$ of
matrices of the form (\ref{mat}) with  $k_1,k_2\in\ZZ$. The first
possibility in fact cannot occur,  since in that case the action
of $U_n$ on $O(p)$, and therefore on $M$, is not effective.

Assume
that $K:=k_1(n-1)-k_2\ne \pm 1,0$. Since $(k_1,k_2)=1$,
 either $k_1$ or $k_2$ is not a multiple of $K$. We set
 $t:=2\pi k_1/K$ in the first case and  $t:=2\pi k_2/K$
in the second case. Then
$e^{it}\cdot\hbox{id}$ is a nontrivial central element of $U_n$ that
belongs
to $H_{k_1,k_2}$. Hence the action of $U_n$ on $O(p)$ is not
effective, which is a contradiction. Further,  assuming
that $K=0$ we obtain $k_1=\pm 1$ and $k_2=\pm(n-1)$. But the center of $U_n$  clearly
 lies  in $H_{1,n-1}$, which yields that the action is
not effective again. Hence $K=\pm 1$.

Now let  $K=-1$. It is not difficult to show that each
 element of the corresponding group
$H_{k_1,k_1(n-1)+1}$ can be expressed in the following form:
\begin{equation}
\left(\begin{array}{cc}
(\det B)^k & 0\\
0 & (\det B)^kB
\end{array}
\right),\label{matf}
\end{equation}
where $B\in U_{n-1}$ and $k:=k_1$. In a similar way, if $K=1$, then each
element of the corresponding group $H_{k_1,k_1(n-1)-1}$ can be expressed
in the form (\ref{matf}) with $k:=-k_1$.

Let $m:=|nk+1|$ and consider the lense manifold ${\cal L}^{2n-1}_m$. We
claim that $O(p)$ is isomorphic to ${\cal
  L}_m^{2n-1}$. We identify $\ZZ_m$ with the subgroup of $U_n$  consisting
of the matrices  $\sigma\cdot\hbox{id}$ with $\sigma^m=1$ and
consider the standard action of $U_n/\ZZ_m$ on ${\cal L}^{2n-1}_m$. The
isotropy subgroup $S$ of the point in ${\cal L}^{2n-1}_m$ represented by
the point
$(1,0,\dots,0)\in S^{2n-1}$ is the standard embedding of $U_{n-1}$ in
$U_n/\ZZ_m$, namely, it consists of elements $C\ZZ_m$, where
$$
C=\left(\begin{array}{cc}
1& 0\\
0& B
\end{array}\right)
$$
and  $B\in U_{n-1}$. The manifold $(U_n/\ZZ_m)/S$ is equivariantly
diffeomorphic to ${\cal L}^{2n-1}_m$. We  now show that it is also
isomorphic to $O(p)$. Indeed, consider the Lie group
isomorphism
\begin{equation}
\phi_{n,m}:U_n/\ZZ_m\ra U_n,\qquad \phi_{n,m}(A\ZZ_m)=(\hbox{det}\,A)^k\cdot A,\label{isom}
\end{equation}
where $A\in U_n$. Clearly, $\phi_{n,m}(S)\subset U_n$ is the subgroup
 of matrices of the form
 (\ref{matf}), that is,  $H_{k_1,k_2}$.
Thus, it is conjugate in $U_n$ to $I_p$, and therefore
$(U_n/\ZZ_m)/S$ is isomorphic to $U_n/I_p$ and
to $O(p)$. More precisely, the isomorphism $f:{\cal L}^{2n-1}_m\ra O(p)$ is
the following composition of maps:
\begin{equation}
f=f_1\circ \phi_{n,m}^*\circ f_2,\label{f}
\end{equation}
where $f_1:U_n/H_{k_1,k_2}\ra O(p)$ and $f_2:{\cal L}^{2n-1}_m\ra (U_n/\ZZ_m)/S$
are the standard equivariant equivalences and the isomorphism
$\phi_{n,m}^*:(U_n/\ZZ_m)/S\ra U_n/H_{k_1,k_2}$ is induced by $\phi_{n,m}$ in
the obvious way. Clearly, $f$ satisfies
\begin{equation}
f(gq)=\phi_{n,m}(g)f(q),\label{equivar1}
\end{equation}
for all $g\in U_n/\ZZ_m$ and $q\in {\cal L}^{2n-1}_m$.

Thus, $f$ is an isomorphism between ${\cal L}^{2n-1}_m$ and $O(p)$ regarded
as
homogeneous spaces, as required. \qed
\smallskip\\

The next result shows that isomorphism (\ref{f}) in Proposition
\ref{realh2} is either a CR or an anti-CR diffeomorphism.

\begin{proposition}\label{realh3}\sl \sl Let $M$ be a
  complex manifold of dimension $n\ge 2$ endowed with  an effective
action of  $U_n$
  by biholomorphic transformations. For
  $p\in M$ suppose that $O(p)$ is a real hypersurface in
  $M$                         isomorphic as a homogeneous space to a
  lense manifold ${\cal L}^{2n-1}_m$. Then an isomorphism ${\cal F}:{\cal
  L}^{2n-1}_m\ra O(p)$ can be chosen to be a CR-diffeomorphism that
  satisfies either the relation
\begin{equation}
{\cal F}(gq)=\phi_{n,m}(g){\cal F}(q),\label{equivar22}
\end{equation}
or the relation
\begin{equation}
{\cal F}(gq)=\phi_{n,m}(\overline{g}){\cal F}(q),\label{equivar33}
\end{equation}
for all $g\in U_n/\ZZ_m$ and $q\in {\cal L}^{2n-1}_m$ (here ${\cal L}^{2n-1}_m$ is
considered with the
  CR-structure inherited from $S^{2n-1}$).
\end{proposition}

\noindent {\bf Proof:} Consider the standard covering map
$\pi:S^{2n-1}\ra{\cal L}^{2n-1}_m$ and  the induced map $\tilde\pi:=f\circ
\pi:S^{2n-1}\ra O(p)$, where $f$ is defined in (\ref{f}).
It follows from (\ref{equivar1})
that the covering map $\tilde\pi$ satisfies
\begin{equation}
\tilde\pi(gq)=\tilde\phi_{n,m}(g)\tilde\pi(q),\label{equivar2}
\end{equation}
for all $g\in U_n$ and $q\in S^{2n-1}$ where
$\tilde\phi_{n,m}:=\phi_{n,m}\circ\rho_{n,m}$ and  $\rho_{n,m}:U_n\ra U_n/\ZZ_m$
is  the
standard projection.

Using $\tilde\pi$  we can pull back the CR-structure from $O(p)$ to
$S^{2n-1}$. We denote by $\tilde S^{2n-1}$ the sphere $S^{2n-1}$
equipped  with
this new CR-structure. It follows from (\ref{equivar2}) that the CR-structure on
 $\tilde
S^{2n-1}$ is invariant under the standard action of $U_n$ on
$S^{2n-1}$.

We  now prove the following lemma.

\begin{lemma}\label{stand}\sl There exist exactly two CR-structures on
  $S^{2n-1}$ invariant under the standard action of $U_n$, namely, the
  standard CR-structure on $S^{2n-1}$ and the structure obtained by
  conjugating the standard one.
\end{lemma}

\noindent {\bf Proof of Lemma \ref{stand}:} For
  $q_0:=(1,0,\dots,0)\in S^{2n-1}$  let $I_{q_0}$ be the isotropy
  subgroup of this point with respect to the standard action of $U_n$ on
  $S^{2n-1}$. Clearly, $I_{q_0}=U_{n-1}$, where $U_{n-1}$ is embedded
  in $U_n$ in the standard way. Let $L_{q_0}$ be the corresponding
  linear isotropy subgroup. Clearly, the only $(2n-2)$-dimensional subspace of
  $T_{q_0}(S^{2n-1})$ invariant under the action of $L_{q_0}$ is
  $\{z_1=0\}$. Hence there exists a unique contact
  structure on $S^{2n-1}$ invariant under the standard action of
  $U_n$.

On the other hand there exist exactly two ways to introduce in
$\RR^{2n-2}$  a  $U_{n-1}$-invariant structure of complex linear space:
the standard complex structure and  its conjugation (this  is obvious
for $n=2$, and easy to show for $n\ge 3$, and therefore we shall omit
the proof). Let  $J_q$ be  the operator of complex structure in the corresponding
subspace of $T_q(S^{2n-1})$, $q\in S^{2n-1}$. Since there exist  only two
possibilities for $J_q$, and $J_q$ depends smoothly on $q$, the
lemma follows.\qed
\smallskip\\

Proposition \ref{realh3} easily follows from Lemma \ref{stand}. Indeed, if the
CR-structure of
$\tilde S^{2n-1}$ is identical to that of $S^{2n-1}$, then
we set ${\cal
  F}:=f$. Clearly, ${\cal F}$ is a CR-diffeomorphism and satisfies
 (\ref{equivar22}).
On the other hand, if the CR-structure of $\tilde S^{2n-1}$ is obtained from
the structure
of $S^{2n-1}$ by conjugation, then we set ${\cal
  F}(t):=f(\overline{t})$ for $t\in {\cal L}^{2n-1}_m$. Clearly,
${\cal F}$ is a CR-diffeomorphism and satisfies (\ref{equivar33}).

The proof of the proposition is complete.\qed
\smallskip\\

We  introduce now additional  notation.

\begin{definition}\label{fachopf} \sl Let
  $d\in\CC\setminus\{0\}$, $|d|\ne 1$, let $M_d^n$ be the Hopf
  manifold constructed by identifying $z\in\CC^n\setminus\{0\}$ with
  $d\cdot z$, and let $[z]$ be the equivalence class of $z$. Then
we denote by
$M_d^n/\ZZ_m$, with $m\in\NN$,
the complex manifold obtained from $M_d^n$ by identifying
$[z]$ and $[e^{\frac{2\pi i}{m}}z]$.
\end{definition}

We are now ready to prove the following theorem.

\begin{theorem}\label{finalstep} \sl Let $M$ be a connected complex manifold
of dimension $n\ge 2$ endowed  with an  effective action of  $U_n$
   by biholomorphic transformations. Suppose that
  all orbits of this action are real
  hypersurfaces.  Then there exists $k\in\ZZ$ such that,
for $m=|nk+1|$, $M$ is   biholomorphically equivalent to
  either

\noindent (i) $S_{r,R}^n/\ZZ_m$, where
  $S_{r,R}^n:=\{z\in\CC^n:r<|z|<R\}$, $0\le r<R\le\infty$, is a spherical
  layer, or

\noindent (ii) $M_d^n/\ZZ_m$.

\noindent The biholomorphic equivalence $f$ can be chosen to satisfy either
the relation
\begin{equation}
f(gq)=\phi_{n,m}^{-1}(g)f(q),\label{a}
\end{equation}
or the relation
\begin{equation}
f(gq)=\phi_{n,m}^{-1}(\overline{g})f(q),\label{a1}
\end{equation}
for all $g\in U_n$ and $q\in M$, where $\phi_{n,m}$ is defined in
(\ref{isom}) (here $S_{r,R}^n/\ZZ_m$ and $M_d^n/\ZZ_m$ are equipped  with
the standard actions of $U_n/\ZZ_m$).
\end{theorem}

\noindent {\bf Proof:} Assume first that $M$
is non-compact. Let $p\in M$.
By Propositions \ref{realh2} and \ref{realh3},
for some $m=|nk+1|$, $k\in\ZZ$,
 there exists a CR-diffeomorphism
$f:O(p)\ra{\cal L}^{2n-1}_m$  such that
either (\ref{a}) or (\ref{a1}) holds
for all $q\in O(p)$. Assume  first that (\ref{a}) holds.
The map $f$ extends
  to a biholomorphic map of a neighborhood $U$ of $O(p)$ onto a
  neighborhood of ${\cal L}^{2n-1}_m$ in $(\CC^n\setminus\{0\})/\ZZ_m$. We can
take  $U$ to be a connected union of
  orbits. Then the extended map     satisfies (\ref{a}) on $U$, and
  therefore maps $U$ biholomorphically onto the quotient of a spherical layer
   by the action of $\ZZ_m$.

Let $D$ be a maximal
  domain in $M$ such that there exists a biholomorphic
  map $f$ from $D$ onto the quotient of a spherical layer by the action
  of $\ZZ_m$ that
  satisfies a relation of the form (\ref{a}) for all
  $g\in U_n$ and $q\in D$. As  shown
  above, such a domain  $D$ exists. Assume that $D\ne M$ and
let $x$ be a boundary point of $D$. Consider the orbit $O(x)$.
Extending a map from $O(x)$ into a lense manifold to a
  neighborhood of $O(x)$ as above, we see that
 the orbits of all points  close to
  $x$ have the same type as $O(x)$. Therefore, $O(x)$ is also
  equivalent to ${\cal L}^{2n-1}_m$.
Let $h:O(x)\ra {\cal L}^{2n-1}_m$ be a CR-isomorphism. It satisfies
  either relation (\ref{a}) or relation (\ref{a1})
for all $g\in U_n$ and $q\in O(x)$.

Assume  first that (\ref{a}) holds
  for $h$.
The map $h$ extends to some
neighborhood $V$ of $O(x)$ that we can assume  to be a connected union of
orbits.  The extended map     satisfies (\ref{a}) on $V$.
For  $s\in
  V\cap D$ we  consider the orbit $O(s)$. The maps $f$ and $h$ take
  $O(s)$ into some surfaces $r_1S^{2n-1}/\ZZ_m$ and
  $r_2S^{2n-1}/\ZZ_m$, respectively, where $r_1,r_2>0$.
Hence $F:=h\circ f^{-1}$ maps
  $r_1S^{2n-1}/\ZZ_m$ onto $r_2S^{2n-1}/\ZZ_m$ and satisfies the relation
\begin{equation}
F(ut)=uF(t), \label{c}
\end{equation}
for all $u\in U_n/\ZZ_m$ and $t\in r_1S^{2n-1}/\ZZ_m$. Let
$\pi_1:r_1S^{2n-1}\ra r_1S^{2n-1}/\ZZ_m$ and $\pi_2:r_2S^{2n-1}\ra
r_2S^{2n-1}/\ZZ_m$ be the standard projections. Clearly,
 $F$ can be lifted to a map     between $r_1S^{2n-1}$ and
$r_2S^{2n-1}$, i.e., there exists a CR-isomorphism $G: r_1S^{2n-1}\ra
r_2S^{2n-1}$ such that
\begin{equation}
F\circ\pi_1=\pi_2\circ G. \label{pod}
\end{equation}
We see from (\ref{c}) and (\ref{pod}) that,
 for all $g\in U_n$ and $y\in r_1S^{2n-1}$,
$$
\begin{array}{l}
\pi_2(G(gy))=F(\pi_1(gy))=F(\rho_{n,m}(g)\pi_1(y))=\\
\rho_{n,m}(g)F(\pi_1(y))=\rho_{n,m}(g)\pi_2(G(y))=\pi_2(gG(y)),
\end{array}
$$
where $\rho_{n,m}:U_n\ra U_n/\ZZ_m$ is the standard projection. Since the
fibers of $\pi_2$ are  discrete, this leads to the relation
\begin{equation}
G(gy)=gG(y),\label{d}
\end{equation}
for all $g\in U_n$ and $y\in r_1S^{2n-1}$.

The map $G$ extends to a biholomorphic
map of the corresponding balls $r_1B^n$, $r_2B^n$, and the
extended map     satisfies (\ref{d}) on $r_1B^n$. Setting $y=0$ in
(\ref{d}) we see that
 $G(0)$ is a fixed point of  the standard action of $U_n$ on
$r_2B^n$, and therefore $G(0)=0$. Combined with (\ref{d}) this shows
that
$G=d\cdot\hbox{id}$,
where $d\in \CC\setminus\{0\}$. This
means, in particular, that $F$ is biholomorphic
on $(\CC^n\setminus\{0\})/\ZZ_m$.
Now,
$$
H:=\Biggl\{\begin{array}{l}
F\circ f\qquad\,\hbox{on $D$}\\
h\qquad\qquad\hbox{on $V$}
\end{array}
$$
is a holomorphic map
on $D\cup V$, provided that $D\cap V$ is connected.

We  now claim that we can choose $V$ such that
 $D\cap V$ is connected. We assume that $V$ is
small enough, hence the strictly pseudoconvex  orbit $O(x)$
partitions  $V$ into two pieces. Namely,   $V=V_1\cup
V_2\cup O(x)$, where  $V_1\cap V_2=\emptyset$ and   each
intersection $V_j\cap D$ is
connected. Indeed, there exist  holomorphic coordinates on $D$ in which
$V_j\cap D$ is a union of the quotients of spherical layers
by the action of
$\ZZ_m$. If there are several  such ``factorized''
layers, then there exists a layer with  closure disjoint
from  $O(x)$ and hence $D$ is disconnected, which is
impossible.  Therefore, $V_j\cap D$ is connected and, if $V$ is
sufficiently small, then each $V_j$ is either a subset of $D$ or is disjoint
from
 $D$. If $V_j\subset D$ for $j=1,2$, then $M=D\cup V$ is
compact which contradicts our assumption. Thus, only one set  of
$V_1$, $V_2$ lies  in $D$, and therefore
$D\cap V$ is connected. Hence the map $H$ is well-defined. Clearly,
it satisfies (\ref{a})
for all $g\in U_n$ and $q\in D\cup V$.

We will now show that $H$ is one-to-one on $D\cup
V$. Obviously, $H$ is one-to-one on each of $V$ and $D$. Assume that there exist
points $p_1\in D$ and $p_2\in V$ such that $H(p_1)=H(p_2)$. Since $H$ satisfies
(\ref{a}) for all $g\in U_n$ and $q\in D\cup V$, it follows that
$H(O(p_1))=H(O(p_2))$. Let $\Gamma(\tau)$, $0\le \tau\le 1$ be a continuous path in
$D\cup V$ joining $p_1$ to $p_2$. For each $0\le \tau\le 1$ we set
$\rho(\tau)$ to be the radius of the sphere
corresponding to the lense manifold
$H(O(\Gamma(\tau)))$. Since $\rho$ is continuous and
$\rho(0)=\rho(1)$, there exists a point $0<\tau_0<1$ at which $\rho$
attains either its maximum or its minimum on $[0,1]$.  Then
 $H$ is not
one-to-one in a neighborhood of $O(\Gamma(\tau_0))$, which is a
contradiction.

We have thus constructed a domain  containing $D$ as a proper subset
 that can be mapped onto the quotient of a spherical layer by the action
of $\ZZ_m$ by means of a map satisfying (\ref{a}).
This is a contradiction showing that in fact $D=M$.

Assume now that $h$ satisfies (\ref{a1}) (rather than (\ref{a}))
  for all $g\in U_n$ and $q\in
O(x)$. Then $h$ extends to a neighborhood $V$ of $O(x)$ and satisfies
(\ref{a1}) there. For  a point $s\in
  V\cap D$ we  consider its orbit $O(s)$. The maps     $f$ and $h$ take
  $O(s)$ into some lense manifolds $r_1S^{2n-1}/\ZZ_m$ and
  $r_2S^{2n-1}/\ZZ_m$, respectively, where  $r_1,r_2>0$.
Hence $F:=h\circ f^{-1}$
maps
  $r_1S^{2n-1}/\ZZ_m$ onto $r_2S^{2n-1}/\ZZ_m$ and satisfies the relation
\begin{equation}
F(ut)=\overline{u}F(t), \label{c1}
\end{equation}
for all $u\in U_n/\ZZ_m$ and $t\in r_1S^{2n-1}/\ZZ_m$. As above,
$F$ can be lifted to a map $G$ from $r_1S^{2n-1}$ into
$r_2S^{2n-1}$.
By  (\ref{c1}) and (\ref{pod}),  for all $g\in U_n$
and  $y\in r_1S^{2n-1}$ we obtain
$$
\begin{array}{l}
\pi_2(G(gy))=F(\pi_1(gy))=F(\rho_{n,m}(g)\pi_1(y))=\\
\overline{\rho_{n,m}(g)}F(\pi_1(y))=\rho_{n,m}(\overline{g})\pi_2(G(y))=\pi_2(\overline{g}G(y)).
\end{array}
$$
As above, this shows that
\begin{equation}
G(gy)=\overline{g}G(y),\label{d1}
\end{equation}
for all $g\in U_n$ and $y\in r_1S^{2n-1}$.

The map $G$ extends to a biholomorphic
map between  the corresponding balls $r_1B^n$,  $r_2B^n$, and the
extended map satisfies (\ref{d1}) on $r_1B^n$. By setting $y=0$ in
(\ref{d1}) we see similarly to the  above that
 $G(0)$ is a fixed point of  the standard action of $U_n$ on
$r_1B^n$, and thus $G(0)=0$. Hence $G=d\cdot U$, where
$d\in\CC\setminus\{0\}$ and $U$ is a unitary matrix. This, however,
contradicts (\ref{d1}), and therefore $h$  cannot
satisfy (\ref{a1}) on $O(x)$.

The proof in the case when $f$ satisfies (\ref{a1}) on $O(p)$ is
analogous to the above. In this case we obtain an extension to the whole
 of $M$  satisfying  (\ref{a1}). This completes the proof
in the case of non-compact $M$.

Assume now that $M$ is compact. We consider  a domain $D$ as above and
assume  first that the corresponding map     $f$ satisfies (\ref{a}). Since
$M$ is compact, $D\ne M$. Let $x$ be a boundary point of $D$, and
consider the orbit $O(x)$. We choose a connected neighborhood $V$ of
$O(x)$ as above, and let $V=V_1\cup V_2\cup O(x)$, where $V_1\cap
V_2=\emptyset$ and each $V_j$ is either a subset of $D$ or is disjoint
from $D$. If one domain  of $V_1$, $V_2$ is disjoint from $D$, then,
arguing as above, we arrive at  a contradiction with the maximality of
$D$. Hence $V_j\subset D$, $j=1,2$, and $M=D\cup O(x)$.

We can now extend $f|_{V_1}$ and $f|_{V_2}$ to biholomorphic maps
$f_1$ and $f_2$, respectively, that are defined on $V$, map it onto
spherical layers factorized by the action of $\ZZ_m$, and
satisfy (\ref{a}) on $V$. Then  $f_1$ and $f_2$ map $O(x)$ onto
$r_1S^{2n-1}/\ZZ_m$ and $r_2S^{2n-1}/\ZZ_m$, respectively, for some
$r_1,r_2>0$. Clearly, $r_1\ne r_2$. Hence $F:=f_2\circ
f_1^{-1}$ maps $r_1S^{2n-1}/\ZZ_m$ onto $r_2S^{2n-1}/\ZZ_m$ and satisfies
(\ref{c}). This shows, similarly to the above,
 that $F(<t>_1)=<d\cdot t>_2$ for all $<t>_1\in r_1S^{2n-1}/\ZZ_m$, where
$d\in\CC\setminus\{0\}$ and $<t>_j\in r_jS^{2n-1}/\ZZ_m$ is  the
equivalence class of $t\in r_jS^{2n-1}$, $j=1,2$.
Since $r_1\ne r_2$, it follows  that $|d|\ne
1$. Now, the map
$$
H:=\Biggl\{\begin{array}{l}
f\qquad\,\,\hbox{on $D$}\\
f_1\qquad \hbox{on $O(x)$}
\end{array}
$$
establishes a biholomorphic equivalence between
$M$ and $M_d^n/\ZZ_m$ and satisfies (\ref{a}).

The proof in the case when $f$ satisfies (\ref{a1}) on $D$ is
analogous to the above. In this case we obtain an extension $H$ that
satisfies (\ref{a1}).

The proof of the  theorem is complete.\qed
\smallskip\\

\section{The Case of Complex Hypersurface Orbits}
\setcounter{equation}{0}

We  now discuss   orbits that are  complex hypersurfaces.
 We start with several examples.

\begin{example}\label{complh}\rm

Let $B_R^n$ be the ball of radius
  $0<R\le\infty$ in $\CC^n$ and let  $\widehat{B_R^n}$ be its blow-up
   at the origin, i.e.,
$$
\widehat
{B_R^n}:=\left\{(z,w)\in B_R^n\times\CC\PP^{n-1}:z_iw_j=z_jw_i,\,\,\hbox{for
      all $i,j$}\right\},
$$
where $z=(z_1,\dots,z_n)$ are the standard coordinates in $\CC^n$ and
$w=(w_1:\dots:w_n)$ are the homogeneous coordinates in
$\CC\PP^{n-1}$. We define an action of $U_n$ on $\widehat{B_R^n}$ as
follows. For $(z,w)\in\widehat{B_R^n}$ and $g\in U_n$ we set
$$
g(z,w):=(gz,gw),
$$
where in the right-hand side we use the standard actions of $U_n$ on $\CC^n$
and $\CC\PP^{n-1}$.  The points $(0,w)\in\widehat{B_R^n}$ form an
 orbit $O$, which is a complex hypersurface biholomorphically equivalent to
$\CC\PP^{n-1}$. All other orbits are real hypersurfaces that are the
boundaries of strongly pseudoconvex neighborhoods of $O$.

We fix $m\in\NN$ and denote by $\widehat{B_R^n}/\ZZ_m$ the quotient
 of $\widehat {B_R^n}$ by the equivalence relation $(z,w)\sim e^{\frac{2\pi
    i}{m}}(z,w)$. Let $\{(z,w)\}\in\widehat{B_R^n}/\ZZ_m$ be  the
equivalence class of $(z,w)\in\widehat{B_R^n}$. We now
define in a natural way an action of
$U_n/\ZZ_m$ on $\widehat{B_R^n}/\ZZ_m$: for
$\{(z,w)\}\in\widehat{B_R^n}/\ZZ_m$ and $g\in U_n$ we set
$$
(g\ZZ_m)\{(z,w)\}:=\{g(z,w)\}.
$$
The points $\{(0,z)\}$ form the unique complex hypersurface
orbit $O$, which is biholomorphically equivalent to $\CC\PP^{n-1}$, and
each real hypersurface orbit is the boundary of a strongly
pseudoconvex neighborhood of $O$.

Now let
 $S_{r,\infty}^n=\{z\in\CC^n:|z|>r\}$, $r>0$, be a spherical layer with
infinite outer radius and let $\widetilde{S_{r,\infty}^n}$ be the union of $S_{r,\infty}^n$ and
the hypersurface at infinity  in $\CC\PP^n$, namely,
$$
\widetilde{S_{r,\infty}^n}:=\bigl\{(z_0:z_1:\dots:z_n)\in\CC\PP^n:(z_1,\dots,z_n)\in
S_{r,\infty}^n,\,z_0=0,1\bigr\}.
$$
We shall equip $\widetilde{S_{r,\infty}^n}$ with the standard action of
$U_n$. For $(z_0:z_1:\dots:z_n)\in\widetilde{S_{r,\infty}^n}$ and $g\in U_n$ we set
$$
g(z_0:z_1:\dots:z_n):=(z_0:u_1:\dots:u_n),
$$
where $(u_1,\dots,u_n):=g(z_1,\dots,z_n)$. The points
 $(0:z_1:\dots:z_n)$ at infinity
 form an orbit $O$,
 which is a  complex hypersurface biholomorphically equivalent to
$\CC\PP^{n-1}$. All other orbits are real hypersurfaces that are the
boundaries of strongly pseudoconcave neighborhoods of $O$.

We fix $m\in\NN$ and denote by $\widetilde{S_{r,\infty}^n}/\ZZ_m$ the
quotient  of
$\widetilde{S_{r,\infty}^n}$ by the equivalence relation
$(z_0:z_1:\dots:z_n)\sim e^{\frac{2\pi i}{m}}(z_0:z_1:\dots:z_n)$. Let
$\{(z_0:z_1:\dots:z_n)\}\in \widetilde{S_{r,\infty}^n}/\ZZ_m$ be     the
equivalence class of $(z_0:z_1:\dots:z_n)\in \widetilde{S_{r,\infty}^n}$.
We  consider $\widetilde{S_{r,\infty}^n}/\ZZ_m$ with the standard action
of $U_n/\ZZ_m$, namely, for $\{(z_0:z_1:\dots:z_n)\}\in \widetilde
{S_{r,\infty}^n}/\ZZ_m$ and $g\in U_n$ we set
$$
(g\ZZ_m)\{(z_0:z_1:\dots:z_n)\}:=\{g(z_0:z_1:\dots:z_n)\}.
$$
The points $\{(0:z_1:\dots:z_n)\}$ form a unique complex hypersurface orbit $O$ which is biholomorphically equivalent to
$\CC\PP^{n-1}$, and each real hypersurface orbit is the
boundary of a strongly pseudoconcave neighborhood of $O$.

Finally, let  $\widehat{\CC\PP^n}$ be the blow-up of $\CC\PP^n$ at the
point $(1:0:\dots:0)\in\CC\PP^n$:
$$
\begin{array}{l}
\widehat{\CC\PP^n}:=\Bigl\{\Bigl((z_0:z_1:\dots:z_n),w\Bigr)\in\CC\PP^n\times\CC\PP^{n-1}:z_iw_j=z_jw_i\\
\hbox{for
      all $i,j\ne 0$, $z_0=0,1$}\Bigr\},
\end{array}
$$
where $w=(w_1:\dots:w_n)$ are the homogeneous coordinates in
$\CC\PP^{n-1}$. We define an action of $U_n$ in
$\widehat{\CC\PP^n}$ as follows. For
$\Bigl((z_0:z_1:\dots:z_n),w\Bigr)\in \widehat{\CC\PP^n}$
and $g\in U_n$ we set
$$
g\Bigl((z_0:z_1:\dots:z_n),w\Bigr):=
\Bigl((z_0:u_1:\dots:u_n),gw\Bigr),
$$
where $(u_1,\dots,u_n):=g(z_1,\dots,z_n)$.  This action has
exactly two  orbits that are complex hypersurfaces: the orbit $O_1$
consisting  of  the points $\Bigl((1:0:\dots:0),w\Bigr)$ and the orbit $O_2$
 consisting  of the points $\Bigl((0:z_1:\dots:z_n),w\Bigr)$. Both $O_1$
and $O_2$ are biholomorphically equivalent to $\CC\PP^{n-1}$. The real
hypersurface orbits are  the boundaries of strongly pseudoconvex
neighborhoods of $O_1$ and strongly pseudoconcave neighborhoods of $O_2$.

We fix $m\in\NN$ and denote by $\widehat{\CC\PP^n}/\ZZ_m$ the quotient of
$\widehat{\CC\PP^n}$ by the equivalence relation
$\Bigl((z_0:z_1:\dots:z_n),w\Bigr)\sim e^{\frac{2\pi i}{m}}
\Bigl((z_0:z_1:\dots:z_n),w\Bigr)$. Let
$\Bigl\{\Bigl((z_0:z_1:\dots:z_n),w\Bigr)\Bigr\}\in \widehat{\CC\PP^n}/\ZZ_m$
be the
equivalence class of $\Bigl((z_0:z_1:\dots:z_n),w\Bigr)\in
\widehat{\CC\PP^n}$. We shall
consider $\widehat{\CC\PP^n}/\ZZ_m$ with the standard action
of $U_n/\ZZ_m$,
namely, for $\Bigl\{\Bigl((z_0:z_1:\dots:z_n),w\Bigr)\Bigr\}\in \widehat{\CC\PP^n}/\ZZ_m$ and $g\in U_n$ we set:
$$
(g\ZZ_m)\Bigl\{\Bigl((z_0:z_1:\dots:z_n),w\Bigr)\Bigr\}:=\Bigl\{g\Bigl((z_0:z_1:\dots:z_n),w\Bigr)\Bigr\}.
$$
As above, there exist exactly two orbits
that are  complex hypersurfaces: the orbit
$O_1$
consisting of the points $\Bigl\{\Bigl((1:0:\dots:0),w\Bigr)\Bigr\}$ and the orbit $O_2$
 consisting of the points $\Bigl\{\Bigl((0:z_1:\dots:z_n),w\Bigr)\Bigr\}$. Both $O_1$
and $O_2$ are biholomorphically equivalent to $\CC\PP^{n-1}$. The real
hypersurface orbits are  the boundaries of strongly pseudoconvex
neighborhoods of $O_1$ and strongly pseudoconcave neighborhoods of
$O_2$.
\end{example}

We  show below that the  complex hypersurface orbits  in Example
\ref{complh} are in fact the only ones that can occur.

\begin{proposition}\label{two}\sl Let $M$ be a connected complex
  manifold of dimension $n\ge 2$ endowed with
an effective action of  $U_n$
  by biholomorphic transformations. Suppose that each  orbit is  a
  real or a complex hypersurface in $M$. Then there exist at most two
  complex hypersurface orbits.
\end{proposition}

\noindent {\bf Proof:} We fix a smooth  $U_n$-invariant
distance function $\rho$ on $M$.
   Let $O$ be an orbit that is a complex hypersurface.
   Consider the $\epsilon$-neighborhood of $U_{\epsilon}(O)$ of
  $O$ in $M$:
$$
U_{\epsilon}(O):=\left\{p\in M:\inf_{q\in
    O}\rho(p,q)<\epsilon\right\}.
$$
If $\epsilon$ is sufficiently small, then the boundary of $U_{\epsilon}(O)$,
$$
\partial U_{\epsilon}(O)=\left\{p\in M:\inf_{q\in
    O}\rho(p,q)=\epsilon\right\},
$$
is a smooth connected real hypersurface in $M$. Clearly, $\partial
U_{\epsilon}$ is $U_n$-invariant, and therefore it is a
union of orbits. If $\partial U_{\epsilon}(O)$ contains an
 orbit that is a real hypersurface, then $\partial U_{\epsilon}(O)$ obviously
coincides with that orbit.

Assume  that $\partial U_{\epsilon}(O)$ contains an orbit 
that is a complex hypersurface.
Then $\partial U_{\epsilon}(O)$ is a union
of such  orbits. It follows from the proof of
Proposition \ref{dim} (see Case 1 there) that if an orbit $O(p)$ is a complex
hypersurface, then $I_p$ is isomorphic to $U_1\times
U_{n-1}$. By Lemma 2.1 of \cite{IK},  $I_p$ is  in fact
conjugate to   $U_1\times U_{n-1}$ embedded in $U_n$ in the standard
way. Hence the action of  the center of $U_n$  on
$O(p)$ is trivial. Thus, the center of $U_n$ acts
trivially on
each complex hypersurface orbit and hence on the entire  $\partial
U_{\epsilon}(O)$. Then its action on
 $M$ is also  trivial, which contradicts the assumption of
the effectiveness of the  action of
$U_n$ on $M$.

Hence, if $\epsilon$ is sufficiently small, then
$U_{\epsilon}(O)$  contains no complex hypersurface orbits other than
 $O$ itself, and the boundary of $U_{\epsilon}(O)$ is a real
hypersurface orbit. Let  $\tilde M$ be the manifold obtained by
removing all complex hypersurface orbits from $M$. Since such an
orbit has a neighborhood  containing no other complex
hypersurface orbits, $\tilde M$ is connected. It is also clear that
$\tilde M$ is
non-compact. Hence, by Theorem \ref{finalstep}, $\tilde M$ can be
 mapped onto $S_{r,R}^n/\ZZ_m$, for some $0\le
r<R\le\infty$, by a biholomorphic map $f$
 satisfying either (\ref{a}) or (\ref{a1}). The manifold  $S_{r,R}^n/\ZZ_m$
has two ends at infinity, and therefore the number of removed complex
hypersurfaces is at most two, which completes the proof.\qed
\smallskip\\

We can now prove the following theorem.

\begin{theorem}\label{complh1}\sl Let $M$ be a connected complex
  manifold of dimension $n\ge 2$ endowed with an
effective
action of  $U_n$
  by biholomorphic transformations. Suppose that each  orbit
of this action is either a
  real or complex hypersurface  and          at least
  one orbit is a  complex hypersurface. Then there exists $k\in\ZZ$ such that,
for $m=|nk+1|$, $M$ is biholomorphically equivalent to either

\noindent (i) $\widehat{B_R^n}/\ZZ_m$, $0<R\le\infty$, or

\noindent (ii) $\widetilde{S_{r,\infty}^n}/\ZZ_m$, $0\le r<\infty$, or

\noindent (iii) $\widehat{\CC\PP^n}/\ZZ_m$.

\noindent The biholomorphic equivalence $f$ can be chosen to satisfy
  either (\ref{a}) or (\ref{a1})
for all $g\in U_n$ and $q\in M$.
\end{theorem}

\noindent {\bf Proof:} Assume first that only one orbit $O$ is a complex
hypersurface. Consider $\tilde M:=M\setminus O$. Since
$\tilde M$ is clearly non-compact, by Theorem \ref{finalstep} there
exists $k\in\ZZ$ such that for
$m=|nk+1|$ and some $r$ and $R$, $0\le r<R\le\infty$,
the manifold  $\tilde M$ is biholomorphically equivalent
to $S_{r,R}^n/\ZZ_m$ by means of a map $f$  satisfying  either (\ref{a})
or (\ref{a1}) for all
$g\in U_n$ and $q\in\tilde M$.
We shall  assume   that $f$ satisfies (\ref{a}) because the latter case
can be dealt with in the same way.

Suppose first that $n\ge 3$. We fix $p\in O$ and consider $I_p$.  We denote for
the moment by  $H\subset U_n$
the standard embedding of $U_1\times U_{n-1}$ in $U_n$. As
mentioned in the proof of Proposition \ref{two}, there exists $g\in
U_n$ such that $I_p=g^{-1}Hg$. For an arbitrary  real hypersurface
orbit $O(q)$ we set
$$
N_{p,q}:=\left\{s\in O(q):I_s\subset I_p\right\}.
$$
Since $I_s$ is conjugate in $U_n$ to a  subgroup $H_{k_1,k_2}$,
where $k_1:= k$ and $k_2=k(n-1)+1\ne 0$ (see (\ref{matf}) in
the proof of Proposition \ref{realh2}), it follows that
$$
N_{p,q}=\left\{s\in O(q): I_s=g^{-1}H_{k_1,k_2}g\right\}.
$$
It is easy to show now that if we fix $t\in N_{p,q}$, then $N_{p,q}=\{ht\}$,
where
$$
h=g^{-1}\left(
\begin{array}{cc}
\alpha & 0\\
0 & \hbox{id}
\end{array}
\right)g,\qquad \alpha\in U_1.
$$

Let $N_p$ be the union of the $N_{p,q}$'s over all
real hypersurface orbits $O(q)$.
Also let  $N_p'$  be the set of  points in
$S_{r,R}^n/\ZZ_m$ whose isotropy subgroup with respect to the standard
action of $U_n/\ZZ_m$ is
$\phi_{n,m}^{-1}(g^{-1} H_{k_1,k_2} g)$
(see (\ref{isom}) for the definition of $\phi_{n,m}$). It is easy to
verify  that $N_p'$ is a complex curve in $S_{r,R}^n/\ZZ_m$ biholomorphically
equivalent to either an annulus of modulus $(R/r)^m$ (if $0<r<R<\infty$),
or a punctured disk
(if $r=0$, $R<\infty$ or $r>0$, $R=\infty$), or
$\CC\setminus0$ (if $r=0$ and $R=\infty$).
 Clearly,
 $f^{-1}(N_p')=N_p$, and hence $N_p$ is a complex curve in
$\tilde M$.

Obviously, $N_p$ is
invariant under the action of $I_p$. By Bochner's theorem there exist
local holomorphic coordinates in the neighborhood of $p$ such that the
action of  $I_p$ is linear in these coordinates  and
coincides  with the action of the linear isotropy
subgroup $L_p$ introduced in the proof of Proposition \ref{dim} (upon the
natural identification of the coordinate neighborhood in question and a
neighborhood of the origin in  $T_p(M)$). Recall
 that
$L_p$ has two invariant complex subspaces in $T_p(M)$: $T_p(O)$ and a
one-dimensional subspace, which correspond in our coordinates
to $O$ and some holomorphic curve.  It can be easily seen that $\overline{N_p}$ is
precisely this curve. Hence
$\overline{N_p}$ near $p$ is an analytic disc with center at $p$,
and therefore $N_p'$ cannot in fact be equivalent to an annulus, and
we have either $r=0$ or $R=\infty$.

Assume first that $r=0$ and $R<\infty$. We consider a holomorphic embedding
$\nu: S_{0,R}^n/\ZZ_m\ra \widehat{B_R^n}/\ZZ_m$ defined by the formula
$$
\nu(<z>):=\{(z,w)\},
$$
where $w=(w_1:\dots:w_n)$ is uniquely determined by the conditions
$z_iw_j=z_jw_i$ for all $i,j$, and $<z>\in(\CC^n\setminus\{0\})/\ZZ_m$
is the equivalence class of $z=(z_1,\dots,z_n)\in\CC^n\setminus\{0\}$.
Clearly, $\nu$ is $U_n/\ZZ_m$-equivariant.
Now let  $f_{\nu}:=\nu\circ f$. We
claim  that $f_{\nu}$ extends to $O$ as a
biholomorphic map of $M$ onto $\widehat{B_R^n}/\ZZ_m$.

Let $\hat O$
be the orbit in $\widehat{B_R^n}/\ZZ_m$
that is a
 complex hypersurface and let  $\hat p\in
\hat O$ be the (unique) point such that
its isotropy subgroup $I_{\hat
  p}$ (with respect to the
action of $U_n/\ZZ_m$ on $\widehat{B_R^n}/\ZZ_m$ as described in Example
\ref{complh}) is $\phi_{n,m}^{-1}(I_p)$.  Then $\{\hat p\}\cup \nu(N'_p)$
is a smooth complex curve. We define the extension $F_{\nu}$ of
$f_{\nu}$ by setting $F_{\nu}(p):=\hat p$ for each $p\in O$.

We  must show that $F_{\nu}$ is continuous at each point  $p\in O$. Let
$\{q_j\}$ be a sequence of points in $M$  accumulating  to
$p$. Since all accumulation points of the sequence $\{F_{\nu}(q_j)\}$
lie  in $\hat O$ and $\hat O$ is compact, it suffices  to
show that each  convergent subsequence $\{F_{\nu}(q_{j_k})\}$ of
$\{F_{\nu}(q_j)\}$ converges to $\hat p$. For every $q_{j_k}$
there exists $g_{j_k}\in U_n$ such that
$g_{j_k}^{-1}I_{q_{j_k}}g_{j_k}\subset I_p$, i.e.,
$g_{j_k}^{-1}q_{j_k}\in\overline{N_p}$. We select a convergent
subsequence $\{g_{j_{k_l}}\}$ and denote its limit by $g$. Then
$\{g_{j_{k_l}}^{-1}q_{j_{k_l}}\}$ converges to $g^{-1}p$. Since
$g^{-1}p\in O$ and $g_{j_{k_l}}^{-1}q_{j_{k_l}}\in\overline{N_p}$, it follows
that
$g^{-1}p=p$, i.e, $g\in I_p$. The map $F_{\nu}$ satisfies (\ref{a}) for
all $g\in U_n$ and $q\in M$, hence $F_{\nu}(q_{j_{k_l}})\in
\overline{N_{\phi_{n,m}^{-1}(g_{j_{k_l}})\hat p}}$, where
$N_{\phi_{n,m}^{-1}(g_{j_{k_l}})\hat p}\subset \widehat{B^n_R}/\ZZ_m$
is constructed similarly   to $N_p\subset\tilde M$. Therefore the limit
  of $\{F_{\nu}(q_{j_{k_l}})\}$
(equal to the limit of
  $\{F_{\nu}(q_{j_k})\}$) is  $\hat p$. Hence
  $F_{\nu}$ is continuous, and therefore holomorphic on $M$. It obviously
  maps $M$ biholomorphically onto $\widehat{B_R^n}/\ZZ_m$.

The case when
 $r>0$ and $R=\infty$  can be treated along the same lines, but one
must consider the holomorphic embedding
 $\sigma: S_{r,\infty}^n/\ZZ_m\ra \widetilde{S_{r,\infty}^n}/\ZZ_m$
such that
$$
\sigma(<z>):=\{(1:z_1:\dots:z_n)\},
$$
the map
$f_{\sigma}:=\sigma\circ f$, and prove that
 $f_{\sigma}$ extends to $O$ as a
biholomorphic map     of $M$ onto $\widetilde{S_{r,\infty}^n}/\ZZ_m$.

If $r=0$ and $R=\infty$, then precisely one of $f_{\nu}$ and
$f_{\sigma}$ extends to $O$, and the extension defines a biholomorphic
map from $M$ to either
$\widehat{\CC^n}/\ZZ_m$, or $\widetilde{S^n_{0,\infty}}/\ZZ_m$.

Let now $n=2$. We fix $p\in O$ and consider $I_p$. There 
exists $g\in U_2$ such that $I_p=g^{-1}Hg$. As above, we
introduce the sets $N_{p,q}$, i.e., for an arbitrary  real hypersurface
orbit $O(q)$ we set
$$
N_{p,q}:=\left\{s\in O(q):I_s\subset I_p\right\}.
$$
Since $I_s$ is conjugate in $U_2$ to a  subgroup $H_{k_1,k_2}$,
where $k_1:= k$ and $k_2=k+1\ne 0$, it follows that
$$
N_{p,q}=\left\{s\in O(q): I_s=g^{-1}H_{k_1,k_2}g\right\}\cup
\left\{s\in O(q): I_s=g^{-1}h_0H_{k_1,k_2}h_0g\right\},
$$
where
$$
h_0:=\left(
\begin{array}{cc}
0 & 1\\
1 & 0
\end{array}
\right),
$$
i.e., for $n=2$, $N_{p,q}$ has two connected components. We denote them $N_{p,q}^1$
and $N_{p,q}^2$, respectively.
It is easy to show now that if we fix $t\in N_{p,q}$, then
$N_{p,q}^1=\{ht\}$ and  $N_{p,q}^2=\{g^{-1}h_0ght\}$,
where
$$
h=g^{-1}\left(
\begin{array}{cc}
\alpha & 0\\
0 & 1
\end{array}
\right)g,\qquad \alpha\in U_1.
$$

We now consider the corresponding sets $N_p^1$ and $N_p^2$. The point
$p$ is the accumulation point in $O$ for exactly one of these sets. As above, we obtain that
either $r=0$, or $R=\infty$.  For example, assume that
 $r=0$ and $R<\infty$.
Let $\hat O$ be the orbit in $\widehat{B_R^2}/\ZZ_m$ that is a
complex hypersurface. There are precisely two points in $\hat O$ whose
isotropy subgroups in $U_2/\ZZ_m$ coincide with
$\phi_{2,m}^{-1}(I_p)$. These points $\hat p_1$ and $\hat p_2$ are the
accumulation points in $\hat O$ of $\nu(N_p^{'1})$ and
$\nu(N_p^{'2})$, where $N_p^{'1},N_p^{'2}\subset S^n_{0,R}/\ZZ_m$ are the
sets of points with isotropy subgroups equal to
$\phi^{-1}_{2,m}(g^{-1}H_{k_1,k_2}g)$ and
$\phi^{-1}_{2,m}(g^{-1}h_0H_{k_1,k_2}h_0g)$ respectively. We
then define the extension $F_{\nu}$ of $f_{\nu}$ by setting
$F_{\nu}(p)=\hat p_1$ if $N_p^1$ accumulates to $p$ and
$F_{\nu}(p)=\hat p_2$ if $N_p^2$ accumulates to $p$. The proof of
the continuity of $F_{\nu}$ proceeds as for $n\ge 3$. The arguments in the
cases $r>0$, $R=\infty$ and $r=0$, $R=\infty$ are analogous to the above.

Assume  now that  two orbits
$O_1$ and $O_2$ in $M$ are  complex hypersurfaces. As above,
we consider the manifold $\tilde M$ obtained from
$M$ by removing $O_1$ and $O_2$.
 For some $k\in \ZZ$,  $m=|nk+1|$, and
some $r$ and $R$, $0\le r<R\le\infty$,  it is
 biholomorphically equivalent
to $S_{r,R}^n/\ZZ_m$ by means of  a map $f$  satisfying either (\ref{a})
or  (\ref{a1}).
Arguments very similar to the ones used above show that in this case
 $r=0$, $R=\infty$, and $f_{\tau}:=\tau\circ f$ extends to a biholomorphic map
$M\ra\widehat{\CC\PP^n}/\ZZ_m$. Here $\tau: (\CC^n\setminus\{0\})/\ZZ_m\ra
\widehat{\CC\PP}^n/\ZZ_m$ is a $U_n/\ZZ_m$-equivariant map defined as
$$
\tau(<z>):=\Bigl\{\Bigl((1:z_1:\dots:z_n),w\Bigr)\Bigr\},
$$
where $w=(w_1:\dots:w_n)$ is uniquely determined from the conditions
$z_iw_j=z_jw_i$ for all $i,j$.

The proof is complete. \qed
\smallskip\\

\section{The Homogeneous Case}
\setcounter{equation}{0}

We  consider now the case when the action of $U_n$  on $M$ is transitive.

\begin{example}\label{hopff}\rm Examples of manifolds on which $U_n$
  acts transitively and effectively are the Hopf manifolds $M_d^n$ (see
  Definition \ref{fachopf}). Let
  $\lambda$ be a complex number
 such that $e^{\frac{2\pi(\lambda-i)}{nK}}=d$ for
  some $K\in\ZZ\setminus\{0\}$.
We define an action of $U_n$ on $M_d^n$ as follows. Let $A\in
U_n$. We can represent $A$ in the form $A=e^{it}\cdot B$, where $t\in \RR$
and $B\in SU_n$. Then we set
\begin{equation}
A[z]:=[e^{\lambda t}\cdot Bz].\label{action}
\end{equation}
 Of course, we must verify that this  action is well-defined.
Indeed, the same element $A\in U_n$ can be also represented in
the form $A=e^{i(t+\frac{2\pi k}{n}+2\pi l)}\cdot(e^{-\frac{2\pi ik}{n}}B)$,
$0\le k\le n-1$, $l\in\ZZ$.
Then formula (\ref{action}) yields
$$
A[z]=[e^{\lambda(t+\frac{2\pi k}{n}+2\pi l)}\cdot e^{-\frac{2\pi ik}{n}}Bz]=
[d^{kK+nKl}e^{\lambda t}\cdot Bz]=[e^{\lambda t}\cdot Bz].
$$
It is also clear that (\ref{action}) does not depend on the choice of representative in the class $[z]$.

The action in question is obviously transitive. It is also effective. For
let
$e^{it}\cdot B[z]=[z]$ for some $t\in\RR$, $B\in SU_n$, and all
$z\in\CC^n\setminus\{0\}$. Then, for some $k\in\ZZ$, $B=e^{\frac{2\pi i
    k}{n}}\cdot\hbox{id}$, and
some $s\in\ZZ$ the following holds
$$
e^{\lambda t}\cdot e^{\frac{2\pi i k}{n}}=d^s.
$$
Using the definition of $\lambda$ we obtain
$$
\begin{array}{l}
t=\frac{2\pi s}{nK},\\
e^{\frac{2\pi i k}{n}}=e^{-\frac{2\pi i s}{nK}}
\end{array}.
$$
Hence $e^{it}\cdot B=\hbox{id}$, and thus the action is effective.

The isotropy subgroup of the point $[(1,0,\dots,0)]$ is $G_{K,1}\cdot
SU_{n-1}$, where $SU_{n-1}$ is embedded in $U_n$ in the standard way
and $G_{K,1}$ consists of all matrices of the form
$$
\left(\begin{array}{cc}
1 & 0\\
0 & \beta\cdot\hbox{id}
\end{array}\right),
$$
where $\beta^{(n-1)K}=1$.

Another example is provided by the manifolds $M_d^n/\ZZ_m$ (see Definition
\ref{fachopf}). Let $\{[z]\}\in M_d^n/\ZZ_m$ be  the equivalence class of
$[z]$. We define an action of $U_n$ on $M_d^n/\ZZ_m$ by
the formula  $g\{[z]\}:=\{g[z]\}$
for $g\in U_n$. This action is clearly transitive; it is also
effective if, e.g., $(n,m)=1$ and $(K,m)=1$.

The isotropy subgroup of the point
$\{[(1,0,\dots,0)]\}$ is $G_{K,m}\cdot SU_{n-1}$, where $G_{K,m}$ consists of
all matrices of the form
\begin{equation}
\left(\begin{array}{cc}
\alpha & 0\\
0 & \beta\cdot\hbox{id}
\end{array}\right),\label{components}
\end{equation}
with $\alpha^m=1$ and $\alpha^K\beta^{K(n-1)}=1$. Note that in this case
every orbit of  the induced action of $SU_n$ is
equivariantly diffeomorphic to the lense manifold ${\cal L}^{2n-1}_m$.

One can consider more general actions by choosing $\lambda$ such that
$e^{\frac{2\pi(\lambda-i)}{n}}=d^K$, but not all such actions are effective.
\end{example}

We shall now describe complex manifolds  admitting effective
transitive actions of $U_n$. It turns out that such  a manifold is
always
biholomorphically equivalent to one of the manifolds $M_d^n/\ZZ_m$. To
prove this  we shall  look    at orbits of
the induced action of $SU_n$.  We require  the following
algebraic lemma first.

\begin{lemma}\label{odnor}\sl Let $G$ be a connected closed subgroup of $U_n$
  of dimension $n^2-2n$, $n\ge 2$. Then either

\noindent (i) $G$ is irreducible as a subgroup of $GL_n(\CC)$, or

\noindent (ii) $G$ is conjugate to $SU_{n-1}$ embedded in
$U_n$ in the standard way, or

\noindent (iii) for $n=3$, $G$ is conjugate to $U_1\times U_1\times
U_1$ embedded in $U_3$ in the standard way, or

\noindent (iv) for $n=4$, $G$ is conjugate to $U_2\times U_2$
embedded in $U_4$ in the standard way.
\end{lemma}

\noindent {\bf Proof:} We start as in the proof
  of Lemma \ref{un}. Since $G$ is compact, it is completely
reducible, i.e., $\CC^n$ splits into a sum of $G$-invariant
pairwise orthogonal complex subspaces, $\CC^n=V_1\oplus\dots\oplus V_m$,
such that the restriction $G_j$ of $G$ to every $V_j$ is irreducible. Let
$n_j:=\hbox{dim}_{\CC}V_j$ (hence $n_1+\dots+n_m=n$) and let
$U_{n_j}$ be the  unitary
transformation group  of $V_j$. Clearly, $G_j\subset U_{n_j}$, and therefore
$\hbox{dim}\,G\le n_1^2+\dots+n_m^2$.
On the other hand  $\hbox{dim}\,G=n^2-2n$, which shows that
$m\le 2$ for $n\ne 3$. If $n=3$, then it is also possible that  $m=3$,
 which means
that $G$ is conjugate to $U_1\times U_1\times U_1$ embedded in
$U_3$ in the standard way.

Now let  $m=2$.  Then either there exists a unitary transformation of
$\CC^n$ such that each  element of $G$ has
 in the new coordinates the
form (\ref{mat1}) with  $a\in U_1$ and $B\in U_{n-1}$ or, for $n=4$,
$G$ is conjugate to $U_2 \times U_2$. We note that, in the first case,
 the scalars $a$
and the matrices $B$, that arise from elements of $G$ in (\ref{mat1})
form compact connected subgroups of $U_1$ and $U_{n-1}$ respectively; we
shall denote them by $G_1$ and $G_2$ as above.

If $\hbox{dim}\,G_1=0$, then $G_1=\{1\}$, and
therefore $G_2=SU_{n-1}$.

Assume  that $\hbox{dim}\,G_1=1$, i.e., $G_1=U_1$. Therefore, $n\ge
3$. Then
$(n-1)^2-2\le\hbox{dim}\,G_2\le (n-1)^2-1$. It follows from Lemma 2.1
of \cite{IK} that, for $n\ne 3$, we have $G_2=SU_{n-1}$. For $n=3$ it is
 also possible that $G_2=U_1\times U_1$, and
therefore $G$ is conjugate to $U_1\times U_1\times U_1$ embedded in
$U_3$ in the standard way. Assume  that $G_2=SU_{n-1}$ and
consider the Lie algebra ${\frak g}$ of
  $G$. It consists of all matrices of the form (\ref{mat2}) with  $b$ an arbitrary
   matrix in   ${\frak {su}}_{n-1}$ and $l(b)$  a linear
function  of the matrix elements of $b$ ranging  in $i\RR$. However,
$l(b)$ must vanish on the commutant of ${\frak {su}}_{n-1}$ which is
 ${\frak {su}}_{n-1}$ itself. Consequently,
 $l(b)\equiv 0$, which contradicts our assumption that $G_1=U_1$.

The proof is complete.\qed
\smallskip\\

We can now prove the following proposition.

\begin{proposition}\label{realh4}\sl Let $M$ be a complex manifold of
  dimension $n\ge 2$ endowed with an effective
  transitive  action of  $U_n$ by biholomorphic transformations. Then there
  exists $m\in\NN$, $(n,m)=1$, such that for each
  $p\in M$ the orbit  $\tilde O(p)$ of the
  induced action of $SU_n$ is a real
 hypersurface in $M$ that is $SU_n$-equivariantly
  diffeomorphic to the lense manifold ${\cal L}^{2n-1}_m$
  endowed  with the standard action of $SU_n\subset U_n/\ZZ_m$.
\end{proposition}

\noindent {\bf Proof:} Since $M$ is homogeneous under the action of
  $U_n$, for every $p\in M$ we have
$\hbox{dim}\,I_p=n^2-2n$. We now
  apply Lemma \ref{odnor} to the identity component $I_p^c$.
Clearly, if $I_p^c$ contains the
  center of $U_n$, then the action of $U_n$ on $M$ is not effective, and
  therefore cases (iii) and (iv) cannot occur. We  claim that
  case (i) does not occur either.

Since $M$ is compact, the group $\hbox{Aut}(M)$ of all biholomorphic
automorphisms of $M$ is a complex Lie group. Hence we can extend the action of
$U_n$ to a holomorphic transitive action of
$GL_n(\CC)$ on $M$ (see \cite{H}, pp. 204--207). Let $J_p$ be the
isotropy subgroup of $p$ with
respect to this action. Clearly, $\hbox{dim}_{\CC}J_p=n^2-n$. Consider
the normalizer $N(J_p^c)$ of $J_p^c$ in $GL_n(\CC)$. It is known from
results of Borel-Remmert and Tits (see Theorem 4.2 in \cite{A2}) that
$N(J_p^c)$ is a parabolic subgroup of $GL_n(\CC)$. We note that
$N(J_p^c)\ne GL_n(\CC)$. For  otherwise $J_p^c$ would be a normal
subgroup of $GL_n(\CC)$. But $GL_n(\CC)$  contains no normal
subgroup of dimension $n^2-n$. Indeed, considering the intersection of
such a subgroup with $SL_n(\CC)$, we would obtain a normal
subgroup of $SL_n(\CC)$ of positive dimension thus arriving at a contradiction. 

All parabolic subgroups of $GL_n(\CC)$ are well-known.
 Let $n=n_1+\dots+n_r$, $n_j\ge 1$, and let
$P(n_1,\dots,n_r)$ be the group of all matrices that have
blocks of sizes $n_1,\dots,n_r$ on the diagonal,
 arbitrary entries  above the blocks,
and zeros below. Then an arbitrary   parabolic subgroup of
$GL_n(\CC)$ is  conjugate to some  subgroup $P(n_1,\dots,n_r)$.

Since the normalizer $N(J_p^c)$ does not coincide
with $GL_n(\CC)$,  it is conjugate to a subgroup
$P(n_1,\dots,n_r)$ with $r\ge 2$. Hence there exists a proper
subspace of $\CC^n$ that is invariant under the action of $N(J_p^c)$,
and therefore under the action of $I_p^c$. Thus,  $I_p^c$ cannot be irreducible.

Hence there exists $g\in U_n$ such that $gI_p^cg^{-1}=SU_{n-1}$,
  where $SU_{n-1}$ is embedded in $U_n$ in the
  standard way. Clearly, the element $g$ can be chosen from
  $SU_n$, and hence $I_p^c$ is contained in $SU_n$ and is conjugate in
  $SU_n$ to $SU_{n-1}$.

Consider now  the orbit $\tilde O(p)$ of a point $p\in M$
under the induced action of
$SU_n$, and let $\tilde I_p\subset SU_n$ be  the isotropy subgroup of $p$
with respect to this action. Clearly, $\tilde I_p=I_p\cap SU_n$. Since
$I_p^c$ lies in $SU_n$, it follows that
 $\tilde I_p^c=I_p^c$. In particular, $\hbox{dim}\,\tilde
I_p=n^2-2n$, and therefore $\tilde O(p)$ is a real hypersurface in $M$.

Assume now that $n\ge 3$. We require  the following lemma.

\begin{lemma}\label{compo}\sl Let $G$ be a closed subgroup of $SU_n$,
  $n\ge 3$, such that $G^c=SU_{n-1}$, where $SU_{n-1}$ is embedded in
 $SU_n$ in
  the standard way. Let $m$ be the number of connected components of
  $G$. Then $G=G_{1,m}\cdot SU_{n-1}$, where
the group $G_{1,m}$ is defined in (\ref{components}).
\end{lemma}

\noindent {\bf Proof of Lemma \ref{compo}:} Let $C_1,\dots,C_m$ be
 the connected   components of $G$ with
$C_1=SU_{n-1}$. Clearly, there exist $g_1=\hbox{id},g_2,\dots,g_m$ in $SU_n$
such that
$C_j=g_jSU_{n-1}$, $j=1,\dots,m$. Moreover,  for each  pair of indices $i,j$
there exists  $k$ such that $g_iSU_{n-1}\cdot
g_jSU_{n-1}=g_kSU_{n-1}$, and therefore
\begin{equation}
g_k^{-1}g_iSU_{n-1}g_j=SU_{n-1}.\label{comp}
\end{equation}
Applying (\ref{comp}) to the vector $v:=(1,0,\dots,0)$, which is preserved by
 the standard embedding of $SU_{n-1}$ in $SU_n$,
we obtain
$$
g_k^{-1}g_iSU_{n-1}g_jv=v,
$$
i.e.,
$$
SU_{n-1}g_jv=g_i^{-1}g_kv,
$$
which implies that $g_jv=(\alpha_j,0,\dots,0)$,
$|\alpha_j|=1$, $j=1,\dots,m$. Hence $g_j$ has the form
$$
g_j=\left(\begin{array}{cc}
\alpha_j & 0\\
0 & A_j
\end{array}\right),
$$
where $A_j\in U_{n-1}$ and $\det A_j=1/\alpha_j$. Since  $A_j$ can
be written in the form  $A_j=\beta_j\cdot B_j$ with $B_j\in SU_{n-1}$, we can
assume without loss of generality
 that $A_j=\beta_j\cdot\hbox{id}$. Clearly,                each
matrix
$$
g_j\cdot\left(\begin{array}{cc}
1 & 0\\
0 & \sigma\cdot\hbox{id}
\end{array}\right)
$$
where $j$ is arbitrary and  $\sigma^{n-1}=1$,
 also belongs to $G$. Further, it is clear that
the parameters
$\alpha_j$, $j=1,\dots,m$, are all distinct and form a finite subgroup
of $U_1$, which is therefore the group of $m$th roots of unity.

Thus, $G=G_{1,m}\cdot SU_{n-1}$, as required. \qed
\smallskip\\

It now follows from Lemma \ref{compo} that if  $n\ge 3$, then for
each  $p\in M$,
$\tilde I_p$ is conjugate in $SU_n$ to one of the groups
 $G_{1,m}\cdot SU_{n-1}$ with
 $m\in\NN$. Hence $\tilde O(p)$ is $SU_n$-equivariantly
diffeomorphic to ${\cal L}^{2n-1}_m$. Clearly, the  $SU_n$-action is
effective on $\tilde O(p)$ only if  $(n,m)=1$.
The integer $m$ does not depend on $p$ since all isotropy subgroups $I_p$
are conjugate in $U_n$. This proves  Proposition \ref{realh4} for $n\ge 3$.

Now let  $n=2$. Since $\tilde O(p)$ is a homogeneous real hypersurface,
it is either strongly
pseudoconvex  or Levi-flat. Assume  that
$\tilde O(p)$ is Levi-flat. Then it is foliated by complex curves.
Let ${\frak m}$ be the Lie algebra of all
holomorphic vector fields on $\tilde O(p)$ corresponding to the
automorphisms of $\tilde O(p)$
generated by the action of $SU_2$. Clearly, ${\frak m}$ is isomorphic
to ${\frak {su}}_2$. Let  $M_p$ be  the leaf of the foliation passing
through $p$, and consider the subspace ${\frak l}\subset{\frak m}$ of
vector fields tangent to $M_p$ at $p$. The  vector fields in
${\frak l}$ remain tangent to $M_p$ at each point $q\in M_p$, and
therefore ${\frak
  l}$ is in fact a Lie subalgebra of ${\frak m}$.
However,  $\hbox{dim}\,{\frak l}=2$  and ${\frak {su}}_2$ has no 2-dimensional
subalgebras. Hence $\tilde
O(p)$ must be  strongly pseudoconvex.

Similarly to  the proof of Proposition
\ref{realh}, we can now show  that $\tilde I_p$ is isomorphic to a subgroup of
$U_1$. This means that $\tilde I_p$ is
a finite cyclic group, i.e., $\tilde I_p=\{A^l,0\le l< m\}$ for some
$A\in SU_2$ and $m\in \NN$ such that $A^m=\hbox{id}$. Choosing new
coordinates in which $A$ is in the 
diagonal form, we see that
 $\tilde I_p$ is conjugate in $SU_2$ to the group of matrices
$$
\left(\begin{array}{cc}
\alpha & 0\\
0 & \alpha^{-1}
\end{array}\right),\qquad \alpha^m=1.
$$
Hence $\tilde O(p)$ is $SU_2$-equivariantly diffeomorphic to
the lense manifold ${\cal L}^3_m$. Clearly,  the action of $SU_2$ is
effective on $\tilde O(p)$ only if $m$ is  odd.
The integer $m$ does not depend on $p$ since all isotropy subgroups $I_p$
are conjugate in $U_2$. This proves  Proposition \ref{realh} for  $n=2$ and
completes the proof in general.\qed
\smallskip\\

We can now establish  the following result.

\begin{theorem}\label{hopfclass}\sl Let $M$ be a complex manifold of
  dimension $n\ge 2$ endowed with an   effective
  transitive  action of  $U_n$ by biholomorphic transformations. Then $M$ is
  biholomorphically equivalent to some  manifold $M_d^n/\ZZ_m$, where
  $m\in\NN$ and $(n,m)=1$. The equivalence $f:M\ra M_d^n/\ZZ_m$ can be
  chosen  to satisfy either the relation
\begin{equation}
f(gq)=gf(q),\label{qua}
\end{equation}
or, for $n\ge 3$, the relation
\begin{equation}
f(gq)=\overline{g}f(q),\label{qua1}
\end{equation}
for all $g\in SU_n$ and $q\in M$ (here $M_d^n/\ZZ_m$ is considered
with the standard action of $SU_n$).
\end{theorem}

\noindent {\bf Proof:}  We claim first that $M$ is biholomorphically
equivalent to some  manifold $M_d^n/\ZZ_m$. For a proof we only need to
show that $M$ is diffeomorphic
  to $S^1\times{\cal L}^{2n-1}_m$ for some $m\in\NN$ such that
  $(n,m)=1$. Then
 biholomorphic equivalence will follow from Theorem 3.1 of \cite{A1}.

Choose $m$ provided by Proposition \ref{realh4}. For $p\in M$ we
consider the $SU_n$-orbit $\tilde O(p)$.
Let $t_0:=\min\{t>0: e^{it}p\in\tilde
O(p)\}$. Clearly, $t_0>0$. For each  point $q\in\tilde O(p)$
there exists $B\in SU_n$ such that $q=Bp$. Hence
\begin{equation}
e^{it_0}q=e^{it_0}(Bp)=(e^{it_0}B)p=(Be^{it_0})p=B(e^{it_0}p),\label{perevod}
\end{equation}
and $e^{it_0}\tilde O(p)=\tilde O(p)$. This shows that $M':=\cup_{0\le
  t <t_0}e^{it}\tilde O(p)$ is a closed submanifold of $M$ of
  dimension $n$.  Since $M$ is connected, it follows that $M'=M$.

Let $p_t:=e^{it}p$, $0\le t\le t_0$. We consider a curve $\gamma:[0,t_0]\ra
M$ such that $\gamma(0)=\gamma(t_0)=p$, $\gamma(t)\in\tilde O(p_t)$
for each $t$, and $\gamma([0,t_0])$ is diffeomorphic to $S^1$. We can
assume that $\tilde I_p=G_{1,m}\cdot SU_{n-1}$, which is also
 the isotropy
subgroup,
 with respect to the standard action of
$SU_n$ on ${\cal L}^{2n-1}_m$,
 of the point $q\in{\cal L}^{2n-1}_m$ represented by the point
$(1,0,\dots,0)\in S^{2n-1}$. Further, for each
$0<t<t_0$, there exists $g_t\in SU_n$ such that $\tilde
I_{\gamma(t)}=g_t\tilde I_p g_t^{-1}$. Clearly,  $\tilde
I_{\gamma(t)}$ is the isotropy subgroup of the point $q_t:=g_t q$ in
${\cal L}^{2n-1}_m$. Hence the map
$$
\phi_t(h\gamma(t))=hq_t,
$$
where $h\in SU_n$, maps
the orbit $\tilde O(p_t)$
 diffeomorphically (and $SU_n$-equivariantly) onto
${\cal L}^{2n-1}_m$, $0\le t\le t_0$ (here we
set $g_0:=g_{t_0}:=\hbox{id}$, $q_0:=q_{t_0}:=q$).

We define now a map
 $\Phi:M\ra S^1\times{\cal L}^{2n-1}_m$. For each  $x\in M$ there
exists a unique $0\le t<t_0$,
 such that $x\in\tilde O(p_t)$. We  set
$$
\Phi(x)=\bigl(e^{\frac{2\pi i t}{t_0}},\phi_t(x)\bigr).
$$
It is clear that $g_t$, and therefore $q_t$ can be chosen
so that $\Phi$
is a diffeomorphism. Hence $M$ is biholomorphically equivalent to one of
the manifolds $M_d^n/\ZZ_m$.

Let $F:M\ra M_d^n/\ZZ_m$ be a holomorphic equivalence. Using $F$, the
action of $SU_n$ on $M$ can be pushed to an action of $SU_n$ by biholomorphic
transformations on $M_d^n/\ZZ_m$. The group $\hbox{Aut}(M_d^n/\ZZ_m)$ of all
biholomorphic automorphisms of $M_d^n/\ZZ_m$ is isomorphic to
$Q_{d,m}^n:=(GL_n(\CC)/\{d^k\cdot\hbox{id},\,k\in\ZZ\})/\ZZ_m$ (this
can be seen, for example, by lifting automorphisms of $M_d^n/\ZZ_m$ to
its universal cover $\CC^n\setminus\{0\}$). Each
maximal compact subgroup of this group is conjugate to a subgroup of
the form $(U_n/\ZZ_m)\times K$,
where $U_n/\ZZ_m$ is embedded in $Q_{d,m}^n$ in the
standard way, and $K$ is isomorphic to $S^1$. The action of $SU_n$ on
$M_d^n/\ZZ_m$ induces an embedding
$\tau:SU_n\ra Q_{d,m}^n$.
Since $SU_n$ is compact, there exists $s\in Q_{d,m}^n$ such that $\tau(SU_n)$ is contained in
$s((U_n/\ZZ_m)\times K)s^{-1}$. However,  there exists no nontrivial
homomorphism from $SU_n$ into $S^1$, and therefore  $\tau(SU_n)\subset
s(U_n/\ZZ_m)s^{-1}$. Since $(n,m)=1$, it follows that
$\tau(SU_n)=sSU_ns^{-1}$, where $SU_n$ in the right-hand side
is embedded in $Q_{d,m}^n$ in the standard way.

We now set
$f:=\hat s^{-1}\circ F$, where $\hat s$ is the automorphism of $M_d^n/\ZZ_m$
corresponding  to $s\in Q_{d,m}^n$.  Pushing now the action of $SU_n$
on $M$ to an action of $SU_n$ on $M_d^n/\ZZ_m$ by means of  $f$ in place
 of $F$,  for the corresponding embedding $\tau_s:SU_n\ra
Q_{d,m}^n$ we obtain the equality $\tau_s(SU_n)=SU_n$, where $SU_n$
in the right-hand side is embedded
in $Q_{d,m}^n$ in the standard way. Thus,  there exists an
automorphism $\gamma$ of $SU_n$ such that
$$
f(gq)=\gamma(g)f(q),
$$
for all $g\in SU_n$ and $q\in M$.

Assume  first that $n\ge 3$. Then each  automorphism of
$SU_n$ has either the form
\begin{equation}
g\mapsto h_0gh_0^{-1},\label{autoform}
\end{equation}
or the form
\begin{equation}
g\mapsto h_0\overline{g}h_0^{-1},\label{autoform1}
\end{equation}
for some fixed $h_0\in SU_n$ (see, e.g., \cite{VO}).
If $\gamma$ has the form (\ref{autoform}),
then considering in place of $f$
the map $q\mapsto h_0^{-1}f(q)$ we obtain a biholomorphic map  satisfying
(\ref{qua}). If $\gamma$ has  the form (\ref{autoform1}), then
considering in place  of $f$ the map $q\mapsto h_0^{-1}f(q)$ we obtain a
biholomorphic map  satisfying (\ref{qua1}).

Let $n=2$. Then each  automorphism of $SU_2$ has  the form
(\ref{autoform}) and arguing as above  we obtain a
biholomorphic map          satisfying (\ref{qua}).

The proof is complete.\qed
\smallskip\\

\begin{remark}\label{improve}\rm For $n\ge 3$ Theorem \ref{hopfclass}
  can be proved without referring to the results in \cite{A1}.
 We note first that the $SU_n$-equivariant
  diffeomorphism between ${\cal L}^{2n-1}_m$ and $\tilde O(p)$ constructed in
  Proposition \ref{realh4} is either a CR     or an anti-CR map
(here we consider ${\cal L}^{2n-1}_m$ is
   with the CR-structure inherited from $S^{2n-1}$).
The corresponding  proof  is similar  to the proof  of Proposition
  \ref{realh3}. We must only  replace $U_n$ and $U_n/\ZZ_m$ by
  $SU_n$ and $\phi_{n,m}$ by the identity map. Further  we argue as in the
  second part of the proof of Theorem \ref{finalstep} for
  compact $M$, replacing there $U_n$ by $SU_n$.

\end{remark}

\begin{remark}\label{can}\rm Ideally, one would like the biholomorphic
  equivalence in Theorem \ref{hopfclass} to be
$U_n$-equivariant, rather than  just $SU_n$-equivariant.
However, as Example
\ref{hopff} shows, there is no canonical transitive action of $U_n$
on $M_d^n/\ZZ_m$.
\end{remark}

\section{A Characterization of $\hbox{\bf C}^n$}
\setcounter{equation}{0}

In this section we apply the results obtained above to prove the following theorem.

\begin{theorem}\label{char} \sl Let $M$ be a connected complex manifold of dimension
  $n$. Assume that $\hbox{Aut}(M)$ and
  $\hbox{Aut}(\CC^n)$ are isomorphic as topological groups. Then $M$
  is biholomorphically equivalent to $\CC^n$.
\end{theorem}

\noindent {\bf Proof:} The theorem is trivial for $n=1$, so we assume that $n\ge 2$. Since
$M$ admits an effective action of $U_n$ by biholomorphic
transformations, $M$ is biholomorphically equivalent to one of the
manifolds listed in Remark \ref{kaup}, Theorem \ref{finalstep}, Theorem
\ref{complh1} and Theorem \ref{hopfclass}. The automorphism
groups of the following manifolds are clearly Lie groups: $B^n$, $\CC\PP^n$,
$S_{r,R}^n/\ZZ_m$ for $r>0$ or $R<\infty$, $M_d^n/\ZZ_m$, 
$\widehat{B_R^n}/\ZZ_m$,
$\widetilde{S_{r,\infty}^n}/\ZZ_m$, $\widehat{\CC\PP^n}/\ZZ_m$.
Since $\hbox{Aut}(M)$ is isomorphic to
$\hbox{Aut}(\CC^n)$ and $\hbox{Aut}(\CC^n)$ is not locally compact,
$\hbox{Aut}(M)$ cannot be isomorphic to a Lie group and hence $M$ is
not biholomorphically equivalent to any of the above manifolds. 

Therefore, $M$ is biholomorphically equivalent to either $\CC^n$,
or $\CC^{n*}/\ZZ_m$, where $\CC^{n*}:=\CC^n\setminus\{0\}$ and
$m=|nk+1|$ for some $k\in\ZZ$. We will now show that the groups
$\hbox{Aut}(\CC^n)$ and $\hbox{Aut}(\CC^{n*}/\ZZ_m)$ are not
isomorphic.

Let first $m=1$. The group $\hbox{Aut}(\CC^{n*})$ consists of exactly
those elements of $\hbox{Aut}(\CC^n)$ that fix the origin. Suppose
that $\hbox{Aut}(\CC^n)$ and $\hbox{Aut}(\CC^{n*})$ are isomorphic and
let $\psi:\hbox{Aut}(\CC^n)\ra\hbox{Aut}(\CC^{n*})$ denote an
isomorphism. Clearly, $\psi(U_n)$ induces an action of $U_n$ on
$\CC^{n*}$, and therefore, by our results above, there is
$F\in\hbox{Aut}(\CC^{n*})$ such that for the isomorphism
$\psi_F:\hbox{Aut}(\CC^n)\ra\hbox{Aut}(\CC^{n*})$,
$\psi_F(g):=F\circ\psi(g)\circ F^{-1}$, we have: either $\psi_F(g)=g$, or
$\psi_F(g)=\overline{g}$ for all $g\in U_n$.

Consider
$U_{n-1}$ embedded in $U_n$ in the standard way, and consider its
centralizer $C$ in $\hbox{Aut}(\CC^n)$, i.e.,
$$
C:=\left\{f\in\hbox{Aut}(\CC^n):f\circ g=g\circ f\,\,\hbox{for all
    $g\in U_{n-1}$}\right\}.
$$
It is easy to show that $C$ consists of maps $f=(f_1,\dots,f_n)$
such that
\begin{equation}
\begin{array}{l}
f_1=az_1+b,\\
f'=h(z_1)z',
\end{array}\label{f1}
\end{equation}
where $z':=(z_2,\dots,z_n)$, $f':=(f_2,\dots,f_n)$, $a,b\in\CC$, $a\ne
0$, $h(z_1)$ is a nowhere vanishing entire function. Similarly,
let $C^*$ be the centralizer of $U_{n-1}$ in
$\hbox{Aut}(\CC^{n*})$. It consists of  maps $f=(f_1,\dots,f_n)$
such that
\begin{equation}
\begin{array}{l}
f_1=az_1,\\
f'=h(z_1)z',
\end{array}\label{f2}
\end{equation}
where $a\in\CC$, $a\ne 0$, $h(z_1)$ is entire and nowhere
vanishing. Clearly, $\psi_F(C)=C^*$.

Let $C'$ and $C^{*'}$ denote the commutants of $C$ and $C^*$
respectively. Clearly, $\psi_F(C')=C^{*'}$. It is easy to check that
$C^{*'}$ consists exactly of all maps of the form (\ref{f2}) where
$a=1$ and $h(0)=1$. In particular, $C^{*'}$ is Abelian. We will now
show that $C'$ is not Abelian. Indeed, consider the following
elements of $C$ (see (\ref{f1})):
$$
\begin{array}{l}
f(z_1,z'):=(z_1+1,z'),\\
g(z_1,z'):=(2z_1,z'),\\
u(z_1,z'):=(z_1+1,e^{z_1}z').
\end{array}
$$
We now see that
$$
\begin{array}{l}
F(z_1,z'):=f\circ g\circ f^{-1}\circ g^{-1}=(z_1-1,z'),\\
G(z_1,z'):=u\circ g\circ u^{-1}\circ
g^{-1}=(z_1-1,e^{\frac{z_1-2}{2}}z').
\end{array}
$$
Clearly, $F,G\in C'$, and we have
$$
\begin{array}{l}
F\circ G=(z_1-2,e^{\frac{z_1-2}{2}}z'),\\
G\circ F=(z_1-2,e^{\frac{z_1-3}{2}}z').
\end{array}
$$
Hence $F\circ G\ne G\circ F$, and thus $C'$ is not
Abelian. Therefore, $C'$ and $C^{*'}$ are not isomorphic. This
contradiction shows that $\hbox{Aut}(\CC^n)$ and
$\hbox{Aut}(\CC^{n*})$ are not isomorphic.

Let now $m>1$. For $z\in\CC^{n*}$ denote as before by
$<z>\in\CC^{n*}/\ZZ_m$ its equivalence class. Let
$$
H_m^n:=\left\{f\in\hbox{Aut}(\CC^{n*}):<f(z)>=<f(\tilde z)>,\,\hbox{if
    $<z>=<\tilde z>$}\right\}.
$$
The group $\hbox{Aut}(\CC^{n*}/\ZZ_m)$ is isomorphic in the obvious
way to $H_m^n/\ZZ_m$. Suppose
that $\hbox{Aut}(\CC^n)$ and $\hbox{Aut}(\CC^{n*}/\ZZ_m)$ are isomorphic and
let $\psi:\hbox{Aut}(\CC^n)\ra\hbox{Aut}(\CC^{n*}/\ZZ_m)$ denote an
isomorphism. Clearly, $\psi(U_n)$ induces an action of $U_n$ on
$\CC^{n*}/\ZZ_m$, and therefore there is
$F\in\hbox{Aut}(\CC^{n*}/\ZZ_m)$ such that
for the isomorphism
$\psi_F:\hbox{Aut}(\CC^n)\ra\hbox{Aut}(\CC^{n*})$,
$\psi_F(g):=F\circ\psi(g)\circ F^{-1}$, we have: either
$\psi_F(g)=\phi_{n,m}^{-1}(g)$, or
$\psi_F(g)=\phi_{n,m}^{-1}(\overline{g})$ for all
$g\in U_n$, where we consider $U_n/\ZZ_m$ embedded in $H_m^n/\ZZ_m$.

The rest of the proof proceeds as for the case $m=1$ above with obvious 
modifications. We consider the centralizer
$C_m^*$ of $\phi_{n,m}^{-1}(U_{n-1})=\phi_{n,m}^{-1}(\overline{U_{n-1}})\subset
H_m^n/\ZZ_m$. Clearly, $\psi_F(C)=C_m^*$. Then we find the commutant
$C^{*'}_m$ of $C_m^*$, and we have $\psi_F(C')=C_m^{*'}$. As above, it turns out
that $C_m^{*'}$ is Abelian. Therefore,
$\hbox{Aut}(\CC^n)$ and $\hbox{Aut}(\CC^{n*}/\ZZ_m)$ cannot be
isomorphic.

The proof is complete.\qed

{\obeylines
Centre for Mathematics and Its Applications
The Australian National University
Canberra, ACT 0200
AUSTRALIA
E-mail address: Alexander.Isaev@anu.edu.au
\hbox{ \ \ }
\hbox{ \ \ }
Department of Complex Analysis
Steklov Mathematical Institute
42 Vavilova St.
Moscow 117966
RUSSIA
E-mail address: kruzhil@ns.ras.ru
}

\end{document}